\input amstex
\input xy
\xyoption{all}
\documentstyle{amsppt}
\document
\magnification=1200
\NoBlackBoxes
\nologo
\hoffset1.5cm
\voffset2cm
\vsize15.5cm

\def\K{\bold{K}}

\def\Q{\bold{Q}}
\def\Z{\bold{Z}}
\def\R{\bold{R}}
\def\F{\bold{F}}
\def\P{\bold{P}}

\def\Hom{\roman{Hom}\,}

\bigskip


\bigskip

\centerline{\bf NORI DIAGRAMS}

\medskip

\centerline{\bf  AND PERSISTENT HOMOLOGY}

\bigskip

\centerline{Yuri I.~Manin, Matilde Marcolli} 

\medskip

\centerline{\it ${}^1$Max--Planck--Institut f\"ur Mathematik, Bonn, Germany,}
\smallskip
\centerline{\it ${}^2$California Institute of Technology, Pasadena, USA} 
\centerline{\it University of Toronto, Canada} \centerline{\it Perimeter Institute for Theoretical Physics, Canada}
\vskip1cm

ABSTRACT. Recently, it was found that there is a remarkable intuitive
similarity between studies in theoretical computer science dealing
with large data sets on the one hand, and categorical methods
of topology and geometry in pure mathematics, on the other.
In this article, we treat  the key notion of persistency from
computer science in the algebraic geometric context involving
Nori motivic constructions and related methods.
We also discuss model structures for persistent topology.

\bigskip

\centerline{INTRODUCTION AND SUMMARY}

\medskip

This paper is a contribution to the emerging research field in which 
computational mathematics dealing with large data bases interacts with
topology, homological algebra, and ``brave new algebra'' of homotopy
theory. 

\smallskip

As a remarkable result of such interaction, (various versions of) the notion of
{\it persistent homology} appeared. According to the informative survey presented in [BuSiSc15], 
the general notion of persistence in computational mathematics was informed by the topological data analysis.
For a general introduction and overview of persistent homology and topological data analysis, see
[BoChYv18], [Car09], [EdHar10], [Za05].

\smallskip
Large data bases are first represented by a family of sampled data at various scales.
Then each sample is structured as a topological/algebraic object say, a simplicial 
space or its chain complex, which are interrelated by a nesting relation. Finally, the invariants of these objects
at all scales are compared,  so that those of them that are persistent 
across a sufficiently large 
range of scales become encoded in a {\it  persistence diagram},
or {\it barcode}. An intuitively transparent picture of this kind, leading to a  multidimensional
bar-encoding of derived category of sheaves on a real finite-dimensional
vector space, was developed very recently in [KaSch17], while to our knowledge 
the earliest introduction of
persistence and bar code diagrams, in the topological/homotopical
context, is due to S. Barannikov, [Bar94].

\smallskip

Our starting point was the observation that conversely, some of the important 
and already well formalised technical tools
of algebraic topology and algebraic geometry may be represented
as the product of intuitive search for ``persistent'' properties of topological spaces/
algebraic varieties/schemes ``observed'' at various scales or from various distances.

\smallskip

In this paper, we focus on the formalism of Nori diagrams and Nori motives 
(see [HuM-S17] and [Ar13]) and show that persistence philosophy
presents them in a new light.

\smallskip
In  mathematical community, existence of a rich ramification of persistency ideas
is not as widely known as it deserves. Hence we hope that our input might
be fruitful.
\bigskip

{\bf 0.1. Diagrams in various contexts.}  We work below {\it in a fixed small universe} as it 
was presented in [KaSch06], Ch.~1. We do not mention the universe
explicitly.

\medskip

{\bf 0.1.1. Definitions.} ([HuM-St17], Def. 7.1.1) A {\it diagram} $D$ is a family, consisting
of two disjoint sets $V(D)$ (vertices), $E(D)$ (edges), and a map $\partial :\, E(D)\to V(D)\times V(D)$,
$\partial (e) =( \partial_{out} (e),\partial_{in} (e))$ (orientation of edges). An oriented edge
is sometimes called {\it an arrow.}

\smallskip

{\it Morphism of diagrams $D_1\to D_2$} consists of two maps $V(D_1)\to V(D_2)$, 
$E(D_1)\to E(D_2)$, compatible with orientations.
\smallskip

A diagram {\it with identities} is a diagram $D$ in which for every vertex $v$, exactly one
oriented edge from $v$ to $v$ is given and called the identity edge $id_v$.
Morphism of diagrams with identities must map identities to identities.

\smallskip

For example, each category $\Cal{C}$ defines a diagram with identities  $D(\Cal{C})$  for which
$$
V(D(\Cal{C})):= \roman{Ob}\,\Cal{C},\quad E(D(\Cal{C})):= \roman{Hom}\,\Cal{C},
$$
and $\partial (f:\, X\to Y) := (X,Y).$

\smallskip

Given a diagram $D$ and a category $\Cal{H}$, any morphism of diagrams $D\to D(\Cal{H})$
is called {\it a representation} of $D$. Of course, representations 
(perhaps, satisfying additional compatibility conditions) themselves are objects of a category/
vertices of its diagram etc. This is the universe where various persistence intuitions reside and
constructions of persistence invariants develop.  

\smallskip

We will start here with a brief description of persistence constructions developed in computer science,
and then give a short survey of Nori's persistency.

\bigskip

{\bf 0.2. Thin categories and diagrams.}  Let $(\Cal{S}, \le )$ be a {\it poset} that is, a set $\Cal{S}$ with reflexive, transitive,
and anti--symmetric binary order relation ([KaSch06], Def.~1.1.3.) It defines a diagram $D$
for which $V(D):= \Cal{S}$, $E(D) :=$ the set of all pairs $(s_1, s_2)$ such that $s_1<s_2$,
oriented from $s_1$ to $s_2$.

\medskip

{\bf 0.2.1. Definition.} A category $\Cal{C}$ is called {\it thin}, if each set $\roman{Hom} (X,Y)$ 
consists of $\le 1$ element. 

\smallskip

Clearly, for such a category $\roman{Ob}\,\Cal{C}$ has the canonical structure of a poset:
$X<Y$ iff $X\ne Y$ and $\roman{Hom} (X,Y)$ is non--empty. Conversely, each poset
defines in this way a thin category in which morphisms in $\roman{Hom} (X,Y)$  are 
equivalence classes of {\it oriented paths}
from $X$ to $Y$. Hence, describing a thin category, one may restrict oneself to 
an explicit description of only {\it generating morphisms} and keep in mind
that each diagram in a thin category is automatically commutative.
Basic examples of posets/thin categories used in data mining are natural numbers $\bold{N}$
and real numbers $\bold{R}$.

\smallskip

Let now $\Cal{I}$ be a category. Then the functors $\Cal{C}\to \Cal{I}$ from a fixed
category $\Cal{C}$ to $\Cal{I}$  form objects of a category denoted $\Cal{I}^{\Cal{C}}$,
with natural transformations as morphisms. If  ``the indexing category'' $\Cal{I}$ is thin,
then  ${\Cal{I}}^{\Cal{C}}$ is also thin. More precisely, a natural transformation
$F \to G$ exists if and only if $F(X)\le G(X)$ for all $X\in \roman{Ob}\, \Cal{C}$,
and this last relation makes from $\Cal{I}^{\Cal{C}}$ a poset.

\smallskip

This remark allows one to define a general analog of the semigroup of oriented translations of the
poset $\bold{R}$: $x\mapsto x+a$, for arbitrary thin category $\Cal{I}$. Namely,
it is   the monoid $\bold{Trans}_{\Cal{I}}:= \Cal{I}^{\Cal{I}}$ with respect to the composition.
It acts on any $\Cal{I}^{\Cal{C}}$ by the precomposition. According to [BuSiSc15], p.~1511,
``we can think of  $\bold{Trans}_{\Cal{I}}$ as a sort of `positive cone' in the monoid
of all endomorphism (i.~e., monotone functions) $\Cal{I}\to \Cal{I}$.''
\medskip
{\bf 0.2.2. Example. Spectral Sequences.} Our exposition below is based upon [GeMa03], pp.~200--218.
\smallskip
Let $r\geq 1$ be an integer. We will call {\it the $r$--th page of a spectral sequence} the following thin indexing
category $\Cal{E}_r$:
$$
\roman{Ob}\,\Cal{E}_r :=   \roman{triples}\ (p,q,r), \roman{where}\  p,q\in \bold{Z}.
$$  
   
Besides identities, a system of generating morphisms of  $\Cal{E}_r$  consists of the arrows 
$$
 d_r^{p,q} :\     (p,q,r) \to (p+r, q-r+1, r).
$$
\smallskip

Moreover, {\it the last page} of a spectral sequence is the thin category $\Cal{E}_{\infty}$
 whose objects are all pairs $(p,q)\in \bold{Z}^2$,
and morphisms are generated by the arrows $(p,q)\to (p+1,q-1)$.  
\medskip

Let now $\bold{A}$ be an additive category, and $F: \Cal{E}_r \to \bold{A}$ be any functor,
satisfying the additional condition
$$
F(d_r^{p,q})\circ F(d_r^{p+r,q-r+1})=0 \ \roman{for\ all} \ p,q.
$$
Such functors form a thin subcategory of $\Cal{E}_r^{\bold{A}}$
an object of which may be called {\it the $r$--th page of an $\bold{A}$--valued
spectral sequence.}

 \medskip
 
 Similarly, for the last page we consider the $\bold{A}$--valued functors $\Cal{E}_{\infty} \to \bold{A}$
 transforming each morphism $(p,q)\to (p+1,q-1)$ into the embedding
$$
 Fil^{p+1}A^n \to Fil^pA^n \ \roman{for}\  n=p+q
$$ 
 where $Fil^*$ is family of filtrations on each object of a sequence of objects $A^n$, $n\in \bold{Z}$, in $\bold{A}$.
 
\medskip

{\bf 0.3. Role of thin diagrams in persistence constructions.} 
Chronologically early definition of {\it persistence module} was {\it a family of vector spaces $V_s$,}
indexed by $s\in\bold{N}$ or by $s\in\bold{R}$, endowed with {\it a family of morphisms} 
$f_{s,t}:\,V_s\to V_t$ whose properties can be succinctly expressed by the statement that
this family forms a functor from {\it the indexing category $\bold{N}$ or $\bold{R}$} to the
{\it target category of vector spaces.} More generally, one can consider functors with values
in a thin category such as  pages  of a spectral sequence.

\medskip

{\bf 0.3.1. Example. Sublevelset  persistence module.} It is a real valued function (say, piecewise continuous)
on a topological space
$f:\, X\to \bold{R}$ considered as a functor $F\in \bold{Top}^{\bold{R}}$, $F:\,t\mapsto f^{-1}(-\infty, t]\subseteq X$.

\medskip

Sublevelset  persistence homology of $f$ is defined as a postcomposition of $F$ and a homology theory.
One can consider points in the indexing diagram $s\in\bold{N}$ or  $s\in\bold{R}$ at which persistence homology
jumps up or down when we increase $t$, say $t^+$ or $t^-$. The resulting sequence of indexed numbers,
together with some additional information about appearing/vanishing homology spaces,
is called {\it the barcode} of this persistent homology. For more details, see Sec.~1.3 below.

\medskip

We omit here an essential construction
of {\it interleaving distance.} It was analysed in categorical terms in [BuSiSc15].
More precisely, the authors have shown that an interleaving distance can be defined  by comparing the monoid   
$\bold{Trans}_{\Cal{I}}$
 with the monoid   $[0,\infty ]$
   by a sublinear projection    $\omega\,: \bold{Trans}_{\Cal{I}}\to[0,\infty ]$:
or with the monoid   $[0,\infty )$
   by a superlinear family   $\Omega\,: [0,\infty)\to \bold{Trans}_{\Cal{I}}$.
   Moreover, the authors observe that   $\omega$ and $\Omega$
are dual in a precise categorical sense. From this observation, many of their properties follow easily.
Then it becomes  clear what is needed to replace the monoids  $[0,\infty ]$ 
and $[0,\infty )$  in order to obtain other ways of measuring interleavings. Comparing 
$\bold{Trans_P}$ with   $[0,\infty ]^n$ and $[0,\infty )^n$ the authors of
[BuSiSc15] show  that the resulting `vector persistence' is stable.

\bigskip

{\bf 0.4. Linear representations of diagrams and Nori's persistence: basic constructions.}
{\it Start} with  the  following data:

\medskip
a) a diagram $D$;
\smallskip
b) a noetherian commutative ring with unit $R$ and the category of finitely generated $R$--modules
$R$--$Mod$;

\smallskip
c) a representation $T$ of $D$ in $R$--$Mod$.
\medskip

Let $End(T)$ be defined as the ring 
$$ End(T):=\{ (\phi_v) \prod_{v\in V(D)} End_R(T(v)) \,|\, \phi_{\partial_{out}(e)}\circ T (e) = T(e)\circ \phi_{\partial_{in}(e)}, \,\, \forall e\in E(D) \}. $$
An inclusion of diagrams $D_1\subset D_2$ such that $T_1=T_2|_{D_1}$ determines a homomorphism
$End(T_2) \to End(T_1)$, by projecting the product $\prod_{v\in V(D_2)} End_R(T_2(v))$ onto the product 
$\prod_{v\in V(D_1)} End_R(T_1(v))$. 

\medskip
{\it Produce}  from the data above the category $C(D,T)$ defined in the following way:

\medskip

d1) If $D$ is finite, then $C(D,T)$ is the category $End(T)$-{\rm Mod} of finitely generated $R$--modules
equipped with an $R$--linear action of $End (T)$.
\smallskip
d2) If $D$ is infinite, first consider its  all finite subdiagrams $F$. 
\smallskip
For each $F$ construct $C(F, T|_F)$ as in d1). Then apply the following limiting procedure.
{\it Objects of} $C(D,T)$ will be all objects of the categories $C(F,T|_F)$. If $F\subset F^{\prime}$,
then each object $X_F$ of $C(F, T|_F)$ gives an object of $X_{F'}$ of
$C(F^{\prime},T|_{F^{\prime}})$, via the map from $End(T_F)$-{\rm Mod} to $End(T_{F'})$-{\rm Mod}
determined by the morphism $End(T_{F'}) \to End(T_F)$ as above.
{\it Morphisms from $X$ to $Y$} in $C(D,T)$ will be  defined as colimits over $F$
of morphisms from $X_F$ to $Y_F$ with respect to these extensions.
\medskip
{\it The result} is called {\it the diagram category $C(D,T)$}.
It is an $R$--linear abelian category which
is endowed with $R$--linear faithful exact forgetful functor 
$$
f_T:\,C(D,T)\to R-Mod.
$$  

For more details, see [HuM-S17], pp.~140-144.

\medskip

{\bf 0.4.1. Universal properties of diagram categories.}   Any representation $T:\,D\to R-Mod$
can be presented as precomposition of the forgetful functor $f_T$
with  an appropriate representation $\tilde{T}:\,D\to C(D,T)$:
$$
T = f_T \circ \tilde{T} .
$$
with the following universal property:

\smallskip

Given any $R$--linear abelian category $A$ with a representation $F:\,D\to A$
and $R$--linear faithful exact  functor $f:\,A\to R-Mod$
with $T=f\circ F$, it factorises through a faithful exact functor
$L(F):\, C(D,T) \to A$ compatibly with decomposition
$$
T = f_T \circ \tilde{T} .
$$

For proofs, cf.~[HuM-S17], pp. 140-141.

\medskip

{\bf 0.4.2. Persistence.} The functor $L(F)$ is actually unique up to unique isomorphism of exact additive functors
 ([HuM-S17], p. 167). It is this functor, constructed for various 
 diagrams of geometric origin in algebraic geometry/topology/...  that is an embodiment of persistency in our context.
 Below we give a sketch of relevant constructions; their development
 in various geometric environments is the content of Sec.~1 of our paper.
 
\medskip

{\bf  0.5. Nori geometric diagrams.} If we have a ``geometric'' category $\Cal{C}$ of spaces/varieties/schemes,
 possibly  endowed with additional
structures, in which one can define morphisms of closed 
embeddings $Y\hookrightarrow X$ (or $Y\subset X$) and morphisms of complements to closed embeddings
$X\setminus Y \to X$, we can  define the Nori diagram of
{\it effective pairs} $D(\Cal{C})$ in the following way (see [HuM-S17], 
pp.~207-208).

\medskip

a). One vertex of $D(\Cal{C})$ is a triple $(X,Y,i)$ where $Y\hookrightarrow X$ is a closed
embedding, and $i$ is an integer.

\medskip

b). Besides obvious identities, there are edges of two types.

\smallskip
b1). Let $(X,Y)$ and $(X^{\prime}, Y^{\prime})$ be two pairs
of closed embeddings.
Every morphism $f: X\to X^{\prime}$ such that $f(Y)\subset Y^{\prime}$
produces  functoriality edges $f^*$ (or rather $(f^*,i)$) going from $(X^{\prime}, Y^{\prime},i)$
to $(X, Y, i)$.

\smallskip

b2). Let $(Z\subset Y\subset X)$ be a stair of closed embeddings. Then
it defines coboundary edges  $\partial$ from $(Y,Z,i)$ to $(X,Y, i+1)$.

\medskip

{\bf  0.5.1. (Co)homological representations of Nori geometric diagrams.} If we start not just from
the initial category of spaces $\Cal{C}$, but rather from a pair $(\Cal{C}, H)$ where $H$ is a cohomology
theory,  then assuming
reasonable properties of this pair, we can define the respective representation $T_H$ of
$D(\Cal{C})$ that we will call  a {\it (co)homological representation of $D(\Cal{C})$}.

\smallskip

For a survey of  such pairs $(\Cal{C}, H)$ that were studied in the context
of Grothendieck's motives, see [HuM-S17], pp. 31-133. The relevant cohomology
theories include, in particular, singular cohomology, and algebraic and holomorphic de Rham cohomologies.

\smallskip

Below we will consider the basic example of cohomological representations of Nori diagrams
that leads to Nori motives.

\medskip

{\bf 0.6. Effective Nori motives ([HuM--S17],  pp. 207--208.)}  Take as a category $\Cal{C}$,
starting object in sec.~2.11 above, the category of
varieties $X$ defined over a subfield $k\subset \bold{C}$.

\smallskip

We can then define the Nori diagram $D(\Cal{C})$ as above. This diagram will be
denoted $Pairs^{eff}$ from now on.

\smallskip

The category of effective mixed Nori motives is the diagram category $C(Pairs, H^*)$ where
 $H^i(X,\bold{Z})$ is the respective singular
cohomology of the analytic space $X^{an}$ (cf. [HuM--S17], pp.~31-34 and further on).

\smallskip

Define the diagram of effective pairs $Pairs^{eff}$ exactly as in the general case .

\smallskip 

It turns out ([HuM--S17], Proposition 9.1.2. p.~208) that the map
$$
H^*:\, Pairs^{eff} \to \bold{Z}-Mod
$$
sending $(X,Y,i)$ to the relative singular cohomology $H^i(X(\bold{C}), Y(\bold{C}); \bold{Z})$,
 naturally extends  to a representation
of the respective Nori diagram in the category of finitely generated abelian groups $\bold{Z}$--$Mod$.

\bigskip

\centerline{1.~NORI GEOMETRIC DIAGRAMS}

\medskip

We start with a detailed exposition of Nori's construction briefly sketched in 0.5.
We extend it by the additional data $(f,\lambda )$ below  following D.~Arapura's
construction of motivic sheaves [Ar13], but tracing his steps in wider
categories of topological spaces.
\medskip
{\bf 1.1. Definition.} {\it The Persistence Diagram $D$ of an appropriate category of topological
spaces has vertices of the form $(f: X \to \R,Y,i,\lambda)$ where

\smallskip
(i) $j:Y\hookrightarrow X$  is a continuous embedding of topological spaces.

(ii) $f: X\to \R$ is a piecewise continuous map with finitely many ``critical values" $t\in \R$.
Criticality here means that the homotopy types of  $X_s =f^{-1}(-\infty,s]$ for $s<t$ and for $s>t$
in a small neighbourhood of $t$ are different.

(iii) $i\in \Z_+$ is a non-negative integer.

(iv) $\lambda\in \R_+$ is a non-negative real number. 
\smallskip
There are three types of edges in $D$:
\smallskip
 (1) Each continuous map $\phi: X \to X'$ such that the diagrams
$$ \xymatrix{ Y \ar[r]^{\phi|_Y} \ar[d]_j & Y' \ar[d]^{j'} \\ X \ar[r]^\phi & X' } \ \ \ \text{ and } \ \ \ 
\xymatrix{ X \ar[rr]^\phi \ar[dr]^f & & X' \ar[dl]_{f'} \\ & \R & } $$

commute, with $\phi|_Y=\phi\circ j$ the restriction, 
gives the corresponding edge $$\phi_*: (f:X\to \R,Y,i,\lambda) \to (f':X'\to \R,Y',i,\lambda).$$
\smallskip
(2) Each pair of inclusions $Z\subset Y \subset X$, with compatible maps to $\R$, produces 
corresponding edge
$$ \partial: (f:X\to \R,Y,i,\lambda) \to (f|_Y:Y\to \R,Z,i-1,\lambda). $$
\smallskip
(3) For any $\lambda \leq \lambda'$, there is an edge
$$ p_{\lambda,\lambda'}: (f:X\to \R,Y,i,\lambda) \to (\ell_\lambda \circ f: X\to \R,Y,i,\lambda'), $$
where $\ell_\lambda : \R \to \R$ is the shift map $t \mapsto t-\lambda$.}

\medskip

This notion of Persistence Diagram of a geometric category is 
not the same as what is usually called a persistence diagram in
the contemporary persistent topology literature (which is a multiset of increasing pairs of
numbers in $\R_+\cup \{ \infty \}$). As we will see in Lemma 1.2.1 below, 
our Persistence Diagrams are closely related to the $(\R,\leq)$--indexed diagrams 
of [BuSc14], hence it is also related to the usual persistence diagrams, 
as we will explain.  

\medskip

{\bf 1.2. Linear representations of persistence diagrams.}
The notion of $(\R,\leq)$--indexed diagrams in the category 
of finite dimensional real vector spaces $\roman{Vec}$ was
considered in [BuSco14].

\smallskip

Objects of the category $\roman{Vec}^{(\R,\leq)}$ of [BuSco14] are functors 
$F:(\R,\leq)\to \roman{Vec}$ from the thin category $(\R,\leq)$ to the category of finite dimensional real vector 
spaces. Its morphisms are natural transformations of such functors. 
It is shown in Section~4 of [BuSco14] that $\roman{Vec}^{(\R,\leq)}$ is an abelian category.

\smallskip

Below we will construct a representation  $T_\R: D \to   \roman{Vec}^{(\R,\leq)}$. 

\smallskip

Start with the following preliminary notations. First, the inclusions 
$\iota^\lambda: X_t \subset X_{t+\lambda}$ of sublevel sets 
$X_t=f^{-1}(-\infty,t]$ and $X_{t+\lambda}=f^{-1}(-\infty, t+\lambda]$ induce
maps of the relative homology groups
$$ 
\iota_{X,i}^\lambda: H_i (X_t, Y_t; \R) \to H_i (X_{t+\lambda}, Y_{t+\lambda}; \R) 
$$
where $Y_t=f|_Y^{-1}(-\infty,t]$ are the induced sublevel sets on $Y$.

Second, an object  $V=(V_t)$ of $\roman{Vec}^{(\R,\leq)}$ is given
by a thin diagram of vector spaces $V=(V_t)$, $t\in \R$.
\medskip

{\bf 1.2.1. Lemma--Definition.} {\it The following maps define 
a representation $T_\R: D \to \roman{Vec}^{(\R,\leq)}$ of the Persistence Diagram  $D$.

\smallskip

A. On objects:
$$ 
 T_\R (f: X\to \R, Y, i,\lambda)_t := V_t:=\roman{Range}(\iota_{X,i}^\lambda)_t .
 $$
 
B. On edges (using notations from Def.~1.1 above): 
$$
T_\R (\phi_i)_t := the\ map\  \roman{Range}(\iota_{X,i}^\lambda)_t \to \roman{Range}(\iota_{X',i}^\lambda)_t .
$$

Furthermore,
$$
T_\R(Z\subset Y\subset X):= the\ map\  \roman{Range}(\iota_{X,i}^\lambda)_t \to \roman{Range}(\iota_{X',i}^\lambda)_t \       
$$
 induced by  the inclusions of sublevel sets
$$ 
H_k (X_t, Y_t; \R) \to H_k (X_{t+\lambda'-\lambda}, Y_{t+\lambda'-\lambda}; \R) 
$$

And finally, for the third type of edges  we have morphisms in homology induced by the inclusions
of sublevel sets
$$ 
H_k (X_t, Y_t; \R) \to H_k (X_{t+\lambda'-\lambda}, Y_{t+\lambda'-\lambda}; \R) .
$$
and the corresponding maps
$$ 
p_{\lambda,\lambda'}: \roman{Range}(\iota_{X,i}^\lambda)_t \to \roman{Range}(\iota_{X,i}^{\lambda'})_{t-\lambda} .
$$
}

This Definition, motivated to a large degree by the algebraic--geometric constructions of [Ar13],
agrees also with the one in [BuSiSc15], Sec.`2.2.4, where the spaces $(V_t)$ above appear as
the persistent homology of $(X,Y)$.

\smallskip

The following remark invokes the main example of persistent homology in the form usually
applied to topological data analysis, \cite{Car09}. Recall that, for a finite set of points
$P$ embedded in a Euclidean space $\R^M$ (or in a more general metric space) the
Vietoris--Rips simplicial complex $K(P)_t$ at scale $t>0$ has $P$ as $0$-skeleton
and has a $k$-simplex for each $(k+1)$-tuple of points $\{ p_0,\ldots, p_k \}\subset P$
such that $\roman{dist}(p,p')\leq t$ for all pairs $p,p'\in \{ p_0,\ldots, p_k \}$. 
\medskip

{\bf 1.2.2. Example}. Let $P\subset \R^M$ be a finite set of points embedded in a Euclidean space (a dataset in
some high dimensional ambient space). An $(\R,\leq)$-diagram of topological spaces (simplicial sets),
that is, a functor $P: (\R,\leq) \to \roman{Top}$, is obtained by taking $P(t)$ to be the Vietoris--Rips 
simplicial complex $K(P)_t$ at scale $t\in \R_+^*$ (and empty for $t<0$). In this case the 
persistent homology as defined above recovers the usual notion of persistent homology of
datasets. 

\medskip

{\bf 1.3. Barcode diagrams.} Now we consider the thin indexing category 
 with objects $n\in\Z$ and morphisms
$\roman{Mor}_{(\Z,\leq)}(n,m)$ consisting of a single morphism for $n\leq m$ and empty
othherwise. We pass to the   category $\roman{Vec}^{(\Z,\leq)}$,  
with objects that are functors $F:  (\Z,\leq)\to \roman{Vec}$ from the category $(\Z,\leq)$ to
the category of finite dimensional real vector spaces and with 
morphisms given by natural
transformations of these functors. The category $\roman{Vec}^{(\Z,\leq)}$ is equivalent to the category of
modules over the ring $\R[x]$ (see Lemma~4.5 of [BuSc14]). 

\smallskip

A finite type object in $\roman{Vec}^{(\R,\leq)}$ is a functor $F: (\R,\leq) \to \roman{Vec}$ such
that $F=\oplus_{j=1}^N \chi_{\Cal{I}_j}$ where $\chi_{\Cal{I}_j}(t)=\R$ for $t\in \Cal{I}_j$ and zero otherwise
and $\chi_{\Cal{I}_j}(t\leq t') =\roman{id}_\R$ if $a,b\in\Cal{I}_j$ and zero otherwise
(Definition~4.1 of [BuSc14]).  
By property (ii) of the functions $f: X \to \R$ in Definition 1.1, 
the sublevel sets $X_t =f^{-1}(-\infty,t]$ have locally constant homotopy type,
so in particular the homology $H_*(X_t,Y_t;\R)$ is locally constant in $t\in \R$.
This implies that the image $T(D)$ in $\roman{Vec}^{(\R,\leq)}$
under the representation of Lemma 1.2.1 consists of finite type objects. In particular,
by Theorem~4.6 of [BuSc14] finite type objects in $\roman{Vec}^{(\R,\leq)}$ are also tame,
that is, all but finitely many values $t\in \Cal{I}$ are regular values, for which there is 
an open interval $\Cal{I} \ni t$ such that $V_t=F(t)$ is constant on $\Cal{I}$. The finitely many
points $t\in \R$ that are not regular values are called critical points. 

\smallskip

The barcode diagram of a finite type object $F$ in $\roman{Vec}^{(\R,\leq)}$ is given
by the multiset of pairs $\{ (a_j,b_j) \}_{j=1}^N$ with $a_j,b_j\in \R\cup \{ \pm \infty \}$
such that $a_j\leq b_j$ and $\{ a_j,b_j\}=\partial \Cal{I}_j$ for $F=\oplus_j \chi_{\Cal{I}_j}$.
The finite $a_j,b_j$ are also the critical points of the object $F$. Let $-\infty < c_0 < c_1< \cdots <
c_M < \infty$ denote the ordered sequence of these critical points.

\medskip

{\bf 1.3.1. Lemma.} {\it The representation $T: \,D\to \roman{Vec}^{(\R,\leq)}$ of Lemma 1.2.1.
determines a representation $T_\Z: D \to \roman{Vec}^{(\Z,\leq)}$. Conversely, the datum
of $T_\Z$ together with the map that assigns to each vertex of $D$ the barcode
diagram of its persistent homology completely determine the representation
$T:D \to \roman{Vec}^{(\R,\leq)}$.}

\smallskip
{\bf Proof.} Let $F_\Z : (\Z,\leq) \to \roman{Vec}$ be the functor that assigns to $n\in \Z$ the vector space $F_\Z(n)$
given by $F(t)$ for $t\in (c_n,c_{n+1})$ for $n=0,\ldots,M-1$, the vector space $F(t)$ for $t>c_M$
for all $n\geq M$ and the vector space $F(t)$ for $t<c_0$ for all $n<0$. To $n\leq m$ the functor
$F_\Z$ assigns the same morphism $F(n\leq m)$. This determines a finite type object $F_\Z$ in the
category $\roman{Vec}^{(\Z,\leq)}$ associated to the finite type object $F$ in $\roman{Vec}^{(\R,\leq)}$. 
It is clear then that knowing $F_\Z$ together with the multiset of points $\{ (a_j,b_j) \}_{j=1}^N$
(the barcode diagram) uniquely determine $F$.

\bigskip

 {\bf 1.4. Diagram Category.}
The representation $T_\Z: D \to \roman{Vec}^{(\Z,\leq)}$ is a representation of the diagram $D$
in the category of $R$--modules for $R=\R[x]$. Thus, we can apply to this representation the
construction of the Nori Diagram Category, see [HuM-S17], Sec.~7.1.2.
\smallskip

Given a representation $T: D \to R$--Mod of a diagram $D$ into the category of $R$--modules
for a commutative ring $R$, the Nori Diagram Category $\Cal{C}(D,T)$ is defined in the following
way
(see 0.4 above and  Sec.~7.1.2 of [HuM--S17]). It is the 
category $\roman{End}(T)$-Mod of modules over the ring of endomorphisms
$$
 \roman{End}(T)=\left\{ (\phi_v)_{v\in V(D)} \,|\,  \phi_v \in \roman{End}_R( T(v) ) \,\ such\ that\, 
\phi_{t(e)} \circ T(e) = T(e) \circ \phi_{s(e)}, \right. 
$$
$$ 
\left. \forall e\in E(D), \,\ with\ source\ and\ target \, s(e), t(e) \in V(D) \right\}. 
$$
The category $\Cal{C}(D,T)$ is an $R$--linear abelian category. 
We denote by $\Phi: \Cal{C}(D,T) \to R$--Mod the forgetful functor. 

\medskip

{\bf 1.4.1. Lemma.} {\it By identifying as above the category $\roman{Vec}^{(\Z,\leq)}$ with $\R[x]$--Mod, we obtain
the Nori Diagram Category $\Cal{C}(D,T_\Z)$ associated to the representation
$T_\Z: D \to \roman{Vec}^{(\Z,\leq)}$ of Lemma 1.3.1, with $\Cal{C}(D,T_\Z)=\roman{End}(T_\Z)$--Mod.}

\bigskip

{\bf 1.5. Persistent homology on the Nori Diagram Category.}
We show here that the persistent homology, constructed as in Lemma 1.2.1,
determines a faithful exact functor $\Cal{C}(D,T_\Z) \to \roman{Vec}^{(\R,\leq)}$ on the
Nori Diagram Category. 

\smallskip

{\bf 1.5.1. Lemma.} {\it
Let $\roman{Vec}_f^{(\R,\leq)}$ denote the full subcategory of $\roman{Vec}^{(\R,\leq)}$
with objects that are of finite type. Then

(1) $\roman{Vec}_f^{(\R,\leq)}$ is an abelian subcategory of $\roman{Vec}^{(\R,\leq)}$.

(2) There is an $\R[x]$--linear faithful exact functor $\Psi: \roman{Vec}_f^{(\R,\leq)} \to \roman{Vec}^{(\Z,\leq)}$
constructed as in Lemma 1.3.1.}

\smallskip

{\bf Proof.}   We have to show that $\roman{Vec}_f^{(\R,\leq)}$ is itself an abelian category and the inclusion functor 
is exact. This is equivalent to showing that $\roman{Vec}_f^{(\R,\leq)}$ is closed under taking kernels 
and cokernels. Let $\alpha: F \to F'$ be a natural transformation of functors $F,F': (\R,\leq) \to \roman{Vec}$
that are of the form $F=\oplus_{k=1}^N \chi_{\Cal{I}_k}$ and $F'=\oplus_{j=1}^M \chi_{\Cal{I}'_j}$. On objects
$t\in \R$ the transformation $\alpha$ acts as an $\R$--linear map $\alpha_t : \oplus_k \chi_{\Cal{I}_k}(t)
\to \oplus_j \chi_{\Cal{I}'_j}(t)$. Then there is a finite collection of points $c_0< \cdots < c_m$ in $\R$,
given by the union of the critical points of $F$ and $F'$ such that, for all $t\in (c_i,c_{i+1})$, with
$c_{-1}=-\infty$ and $c_{m+1}=+\infty$, the map $\alpha_t$ is a linear map 
$\alpha_t: \R^{N_i} \to \R_{M_i}$. Since on morphisms the functor $F$ is determined by
$\chi_{\Cal{I}_k}(t\leq t')=\roman{id}$ if $t,t'\in \Cal{I}_k$ and zero otherwise, and similarly for $F'$, 
the natural transformation diagrams
$$
 \xymatrix{ F(t) \ar[r]^{F(t\leq s)} \ar[d]^{\alpha_t} & F(s) \ar[d]^{\alpha_s} \\
F'(t) \ar[r]^{F'(t\leq s)} & F'(s) } 
$$
imply that $\alpha_t$ is locally constant. Thus, the kernel and cokernel of $\alpha_t$
also determine a finite type object in $\roman{Vec}_f^{(\R,\leq)}$.

\smallskip

Now, consider a  functor $F: (\R,\leq) \to \roman{Vec}$ which is a finite type object in $\roman{Vec}_f^{(\R,\leq)}$.
Proceeding as above, denote by $r_0 < \cdots < r_\ell$  the critical points of $F$, so that $F(t)$ is
locally constant with $F(t)=\R^{N_i}$ for all $t \in (r_i, r_{i+1})$, 
with $r_{-1}=-\infty$ and $r_{\ell+1} = +\infty$. We can then assign to $F$ a functor
$F_\Z: (\Z,\leq) \to \roman{Vec}$ with $F_\Z(n)=F(t)$ with $t\in (r_n,r_{n+1})$ for $n=0,\ldots, \ell-1$
and $F_\Z(n)=F(t)$ with $t<r_0$ for all $n<0$ and $F_\Z(n)=F(t)$ with $t>r_\ell$ for $n\geq \ell$.
We also define on morphisms $F_\Z (n\leq m) = F(t\leq t')$ for $t\in (r_n,r_{n+1})$ (or respectively $t< r_0$ or $t>r_\ell$ depending on the value of $n$) and $t'\in (r_m,r_{m+1})$ (or respectively $t'< r_0$ or $t'>r_\ell$ depending on the value of $m$). Let $\alpha: F \to F'$ be a natural transformation of functors 
$F,F': (\R,\leq) \to \roman{Vec}$. We obtain a corresponding natural transformation $\alpha_\Z: 
F_\Z \to F'_\Z$ by assigning $\alpha_{\Z,n}:  F_\Z(n)\to F'_\Z(n)$ to be the same map
$\alpha_t: F(t) \to F'(t)$ for $t\in (r_n,r_{n+1})$ (or $t< r_0$ or $t>r_\ell$ depending on $n$). 
The transformation $\alpha_{\Z,n}:  F_\Z(n)\to F'_\Z(n)$ is trivial if and only if 
 $\alpha_t: F(t) \to F'(t)$ for $t$ in the corresponding interval
 is also trivial, 
 so that $\alpha_\Z$ is trivial iff $\alpha$ is,
 hence the functor $F\mapsto F_\Z$ is faithful. 
 
 \smallskip
 
 Let $\tilde{\Cal{I}}_n$ denote the intervals $(-\infty,r_0)$ for $n<0$, $\tilde{\Cal{I}}_n=(r_n,r_{n+1})$ for
 $n=0,\ldots,\ell-1$ and $(r_\ell,\infty)$ for $n\geq \ell$.  We have an $\R[x]$-linear structure on 
 $\roman{Vec}_f^{(\R,\leq)}$ where $x$ acts on $F(t)$ as $F(t\leq t')$ for $t\in \tilde{\Cal{I}}_n$ and
 and any $t'\in \tilde{\Cal{I}}_{n+1}$. With respect to this $\R[x]$-linear structure the functor
 $F \mapsto F_\Z$ is $\R[x]$-linear. Moreover, by an argument similar to the one
 used above to check faithfulness, if we have an exact sequence
$$
\xymatrix{0 \ar[r]  & F \ar[r]^\alpha & F' \ar[r]^\beta  & F'' \ar[r] &0}
$$
 in $\roman{Vec}_f^{(\R,\leq)}$ we also obtain a corresponding exact sequence
  $$
  \xymatrix{0 \ar[r] & F_\Z \ar[r]^{\alpha_\Z} & F'_\Z \ar[r]^{\beta_\Z} & F''_\Z \ar[r] & 0}.
  $$
 Hence the functor $\Psi$ mapping $F \mapsto F_\Z$ and $\alpha\mapsto \alpha_\Z$ is
 a faithful exact $\R[x]$--linear functor $\roman{Vec}_f^{(\R,\leq)} \to \roman{Vec}^{(\Z,\leq)}$.

\bigskip

{\bf 1.5.2. Proposition.}  {\it Persistent homology determines a faithful exact functor 
$$
PH_* : \Cal{C}(D,T_\Z) \to \roman{Vec}_f^{(\R,\leq)}.
$$
}
{\bf Proof.}
The Nori Diagram Category $\Cal{C}(D,T)$ of a representation $T: D \to R$--Mod of a diagram $D$
satisfies the following universal property: given any $R$-linear
abelian category $\Cal{A}$, a representation $T_\Cal{A}: D \to \Cal{A}$,  
and an $R$--linear faithful exact functor $\Psi: \Cal{A} \to R$--Mod such that $\Psi\circ T_\Cal{A}=T$,
then there exists a faithful exact functor $\Phi_\Cal{A}: \Cal{C}(D,T)\to \Cal{A}$ such that the following 
diagram commutes:
$$ 
\xymatrix{  & \Cal{C}(D,T) \ar[dd]_(.3){\Phi_\Cal{A}}\ar[dr]^\Phi & \\
D \ar[rr]^{\quad T} \ar[ur]^{\widetilde T} \ar[dr]^{T_\Cal{A}}& & R\roman{-Mod} \\
& \Cal{A} \ar[ur]^\Psi &
} 
$$
As we mentioned in Sec.~0.4 above, the category $\Cal{C}(D,T)$ is in fact completely characterized by this
property up to unique equivalence of categories. (see Sec. 7.1.3 of [HuM-S17], Sec.~7.1.3.)

\medskip

Now,  apply this universal property of the Nori Diagram Category to the following case:
$D$ is the Persistence Diagram of Definition 1.1;
$T_\Z : D \to \R[x]$--Mod is the representation of Lemma 1.5.1; 
$\Cal{A}= \roman{Vec}_f^{(\R,\leq)}$, with the representation $T: D \to \roman{Vec}_f^{(\R,\leq)}$
of Lemma 1.2.1, and the functor $\Psi: \roman{Vec}_f^{(\R,\leq)}\to \R[x]$--Mod is the one
of Lemma 1.5.1. To this purpose it suffices to check that the functors $\Psi, T, T_\Z$
satisfy the composition property $\Psi\circ T =T_\Z$, which is true by construction
(compare Lemma 1.5.1 with Lemma 1.4.1).  We then obtain a
faithful exact functor $PH_* : \Cal{C}(D,T_\Z) \to \roman{Vec}_f^{(\R,\leq)}$ that completes the
commutative diagram 
$$ 
\xymatrix{  & \Cal{C}(D,T_\Z) \ar[dd]_(.3){PH_*} \ar[dr]^\Phi & \\
D \ar[rr]^{\quad T_\Z} \ar[ur] \ar[dr]^{T} & & \R[x]\roman{-Mod} \\
& \roman{Vec}_f^{(\R,\leq)} \ar[ur]^\Psi &
} 
$$
\bigskip

{\bf 1.6. The product structure.}
We will show here that the Persistence Diagram $D$ of Definition 1.1
has the structure of a graded diagram with a commutative product with unit, in
the sense of Definition~8.1.3 of [HuM-S17]. Recall from this definition that
a graded diagram $D$ is a diagram endowed with a map $\deg: V(D) \to \Z/2\Z$
extended to $\deg: E(D)\to \Z/2\Z$ by $\deg(e)=\deg(s(e))-\deg(t(e))$. 
The product $D\times D$  is the diagram with vertices the pairs $(v,w)\in V(D)\times V(D')$
and edges of the form $(e,id)$ or $(id,e')$. A product structure on $D$ is a map of
graded diagrams (a degree preserving map of directed graphs) $D\times D \to D$
together with a choice of edges
$$ 
\alpha_{v,w}: v\times w \to w\times v, \ \ \forall v,w\in V(D) 
$$
$$ 
 \beta_{v,w,u}: v\times (w\times u) \to (v\times w)\times u, 
 $$
 $$
 \beta'_{v,w,u}: (v\times w)\times u \to v\times (w\times u),  . 
 $$
 for all  $v,w,u\in V(D)$.
 A unit is a vertex $\bold{1}$ with $\deg(\bold{1})=0$ and edges $u_v: v \to \bold{1}\times v$ for
 all $v\in V(D)$.

\bigskip

{\bf 1.6.1. Lemma.} {\it
The Persistence Diagram $D$ of Definition 1.1 is a graded diagram 
with  commutative product and unit.}

\smallskip

{\bf Proof.}  Define the $\bold{Z}_2$--grading by $\deg(f:X\to \R,Y,k,\lambda)=k$ mod $2$. 
The product $D\times D \to D$ is given by
$$ 
(f:X\to \R, Y, k,\lambda) \times (f':X'\to \R,Y',k',\lambda') := 
$$
$$ 
(X\times_\R X' \to \R, X \times_\R Y' \cup Y\times_\R X', k+k', \lambda+\lambda'), 
$$
where $X\times_\R X'$ is the fibered product:
$$
 \xymatrix{ X\times_\R X' \ar[r] \ar[d] & X \ar[d]^f \\ X' \ar[r]_{f'} & \R }. 
 $$
 
This product satisfies the identities  
$$
(X\times_\R X')_t=\{ (x,x')\in X\times X'\,:\, f(x)=f'(x')\leq t\} =X_t\times_\R X'_t.
$$
The unit vertex is given by $(\roman{id}: \R\to \R, \emptyset, 0,0)$. The edges $\alpha_{v,w}$,
$ \beta_{v,w,u}$ and $\beta'_{v,w,u}$ are the natural homeomorphisms of topological spaces 
compatible with the maps. 
This finishes the proof.

\medskip

We recall now that in the situation of Lemma 1.6.1 one can define a subclass
of representations 
of $D$ that are compatible with grading and commutative
product. Namely, according to the Definition 8.1.3 of [HuM-S17], the compatibility
 conditions for a representation $T\to R$--Proj of a
graded diagram $D$ are given by the existence of isomorphisms

$$ 
\xymatrix{\tau_{v,w}: T(v\times w) \ar[r]^{\simeq} & T(v)\otimes T(w) }
$$ 
for all $v,w\in V(D)$, with the following properties:
$$
\xymatrix{T(v)\otimes T(w) \ar[r]^{\tau_{v,w}^{-1}} & T(v\times w) \ar[r]^{T(\alpha_{v,w})}  &  T(w\times v) \ar[r]^{\tau_{w,v}} 
& T(w)\otimes T(v)}
$$
is equal to multiplication by $(-1)^{\deg(v) \deg(w)}$; the $\beta$--maps satisfy $T(\beta_{v,w,u})^{-1}=T(\beta'_{v,w,u})$, 
and moreover
$$
\tau_{v,w'}\circ T(1,e) = (\roman{id}\otimes T(e))\circ \tau_{v,w} :\, T(v\times w) \to T(v)\otimes T(w'),
$$
$$
\tau_{v',w}\circ T(e,1) = (T(e)\otimes \roman{id})\circ \tau_{v,w} :\, T(v\times w)\to T(v')\otimes T(w),
$$
$$ \xymatrix{ T(v \times (w\times u)) \ar[r]^{T(\beta_{v,w,u})}
\ar[d]^{\tau\circ \tau} & T((v\times w)\times u) \ar[d]^{\tau\circ\tau} \\ 
T(v)\otimes (T(w)\otimes T(u)) \ar[r]^{\simeq} & (T(v)\otimes T(w))\otimes T(u)) } 
$$
and similarly for the inverse $T(\beta'_{v,w,u})$.

\bigskip
{\bf 1.7. Good persistence vertices.}
In order to define a tensor structure on the Nori Diagram Category $\Cal{C}(D,T_\Z)$
of the Persistence Diagram we need to proceed in a way similar to that adopted in
the construction of the category of Nori motives, see [HuM-S17], Sec.~9. Indeed,
because of the K\"unneth product formula 
$$ 
H_k( (X\times_\R X')_t, (X\times_\R Y'\cup Y\times_\R X')_t;\R) \simeq \oplus_{i+j=k} H_i(X_t, Y_t;\R)\otimes
H_j(X'_t,Y'_t;\R) 
$$
where $(X\times_\R X')_t=X_t\times_\R X'_t$ and $(X\times_\R Y'\cup Y\times_\R X')_t=
X_t\times_\R Y_t'\cup Y_t\times_\R X'_t$,
this relative homology is compatible with the product structure on $D$ in the case where 
these homology groups are supported in a single degree. 
As in the case of Nori motives, we can introduce a class of
``good objects" for which the persistent homology is concentrated in a single degree.

\medskip

{\bf 1.7.1. Definition.} {\it
A vertex $(f:X\to \R, Y, k,\lambda)$ of the Persistence Diagram $D$ is 
a ``good persistence vertex" if the persistent homology 
$$ HP_j(f:X\to \R, Y, k,\lambda)_t= \roman{Range} (H_j(X_t,Y_t;\R) \to H_j(X_{t+\lambda},Y_{t+\lambda};\R)) $$
satisfies $HP_j(f:X\to \R, Y, k,\lambda)_t=0$ for all $j\neq k$.}

\smallskip

Simple examples of good persistence vertices can be constructed as follows.
Let $X$ be a smooth $n$-dimensional compact manifold and let $f: X\to \R$ be a Morse function, which 
has finitely many critical points $x_1,\ldots, x_m$ in $X$ with critical values $c_1<\cdots<c_m$ in $\R$.
The sublevel sets $X_t=f^{-1}(-\infty,t]$ have homotopy type that remains constant when $t$ varies in
each of the intervals $(-\infty,c_1)$, $(c_k,c_{k+1})$ with $k=1,\ldots,m-1$, and $(c_m,\infty)$ and changes
across the critical values by a handle attachment. Let $B$ be an open $n$-ball in $X$ that does not
contain any critical point and such that the sublevel sets $B_t=f|_B^{-1}(-\infty,t]$ are either empty or
a contractible set that is open in the induced topology on $X_t$, or all of $B$. Let $Y=X\smallsetminus B$. 
Then the relative homology $H_k(X_t, Y_t;\R)=H_k(X_t,X_t\smallsetminus B_t;\R)$ is a local
homology and is trivial for $k\neq 0$ and is either trivial or a single copy of $\R$ for $k=n$.
This gives an example where the persistent homology is concentrated in degree $k=n$.

\smallskip

Indeed, if we assume all the topological spaces involved are CW complexes, it
is always possible to compute the homology via a cellular filtration by good spaces.
Indeed we have the following result (Theorem~2.35 of \cite{Hat02}).

\medskip

{\bf 1.7.2. Lemma.} {\it
Given a CW complex $X$ with skeleta $X^{(n)}$, the homology $H_k(X;\Z)$
is computed as the homology of a complex
$$ \cdots \to H_j (X^{(j)}, X^{(j-1)};\Z) \to H_j (X^{(j-1)}, X^{(j-2)};\Z) \to \cdots $$
with the maps given by the boundary maps of the pair of inclusions $X^{(j-2)}\subset X^{(j-1)}\subset X^{(j)}$.
The relative homology $H_n (X^{(j)}, X^{(j-1)};\Z)$ is trivial for $n\neq j$ and a free abelian group for $n=j$
spanned by the $j$-cells of $X$.}

\medskip

{\bf 1.7.3. Lemma.} {\it
Under the assumption that all spaces considered are cellular with cellular maps,
the representations $T: D \to \roman{Vec}^{(\R,\leq)}$ and $T_\Z: D \to \roman{Vec}_f^{(\Z,\leq)}$
of Lemma 1.2.1 and Lemma 1.3.1 are unital graded multiplicative representations.}
\medskip

{\bf Proof.} First notice that the image of $T_\Z: D \to \roman{Vec}^{(\Z,\leq)}=\R[x]$-Mod lies in the
subcategory $\roman{Vec}_f^{(\Z,\leq)}$ of finite type. This is in fact a subcategory of
the category $\R[x]$--Proj of finite projective module over $\R[x]$. Then observe that,
under the cellular assumption, 
Lemma 1.7.2 implies that, if we consider the Nori Diagram Category built on a subdiagram
of the Persistence Diagram where all the vertices are good persistence vertices,
the resulting $C(D^{\roman{good}},T_\Z)$ contains all the objects $(X^{(j)},X^{(j-1)},j,\lambda)$
and $(Y^{(j)}, Y^{(j-1)},j,\lambda)$ for every vertex $(f:X\to \R,Y,j,\lambda)$ of the Persistence Diagram.
Thus, there is an object in the Nori Diagram Category $C(D^{\roman{good}},T_\Z)$ whose image
under the forgetful functor to $\roman{Vec}^{(\Z,\leq)}$ is the same as the image
$T_\Z(f:X\to \R,Y,j,\lambda)$. This implies that we can equivalently use the categories 
$C(D,T_\Z)$ and $C(D^{\roman{good}},T_\Z)$. Using the latter, we can define the product structure,
finishing the proof.

\smallskip

Notice that essentially the same argument was used in [HuM-S17], Sec. 9 in order to construct
the product structure on Nori effective motives. The argument is simplified here
because we work in a topological setting, hence we can directly use 
cellular homology as in Lemma 1.7.2, instead of having to use Beilinson's
fundamental lemma for the cohomology of affine varieties and complexes of
varieties to pass from affine to more general varieties.

\bigskip
{\bf The Tannakian formalism.}
An advantage of reformulating the categorical construction of persistent homology of
[BuSco14], [BuSiSco15], [BluLes17] in terms of the formalism of Nori
diagrams and Nori motives, as we did in the previous subsections, is the fact that this
formulation comes endowed with natural symmetries associated to persistent
homology which are not immediately visible otherwise, namely the associated Tannakian 
formalism. 

\smallskip

In the category of Nori motives, one passes from effective motives to the
localization with respect to $(\bold{G}_m,\{ 1 \},1)$ (inverting the Lefschetz motive) to
obtain a {\it rigid} abelian tensor category to which the Tannakian formalism can be applied.
In our setting we work with weaker properties, as we will discuss more in 
Section 4 where we present a more general formalism
based on Nori diagram for persistent phenomena. We do not assume
that the category we construct is a rigid tensor category, although we will
assume that it has a tensor structure, obtained via approximations using
filtrations by good objects as explained above. Thus, instead of the group
scheme that one expects to obtain as Tannakian Galois group
in the case of rigid tensor categories, we only have a monoid scheme,
obtained as follows.

\smallskip

{\bf 1.8.1. Proposition.}{\it
The representation $T_\Z: D \to \roman{Vec}^{(\Z,\leq)}_f$ induces an equivalence between
the Nori Diagram Category $\Cal{C}(D,T_\Z)$ and the category of finitely generated comodules 
over a bialgebra $\Cal{A}(D,T_\Z)$, which defines a pro--algebraic monoid scheme 
$\roman{Spec}(\Cal{A}(D,T_\Z))$.}

\smallskip

{\bf Proof.} As before, we view the representation $T_\Z: D \to \roman{Vec}^{(\Z,\leq)}_f$
as taking values in the category $\R[x]$-Proj of finitely generated projective modules 
over the Dedekind domain $\R[x]$. We can then apply Theorem~7.1.12 of [HuM-S17]
and we obtain that the Nori Diagram Category $\Cal{C}(D,T_\Z)$ is equivalent to the
category of finitely generated comodules over the coalgebra $\Cal{A}(D,T_\Z)$ given
by the colimit 
$$ 
\Cal{A}(D,T_\Z) = \roman{colim}_{D_F} \roman{End}(T_\Z |_{D_F})^\vee 
$$
over finite sub-diagrams $D_F$
with $\vee$ the $\R[x]$-dual. 

Indeed, as shown in Sec. 7.5.1 of [HuM-S17], if
$R$ is a Dedeking domain, then for the $R$-algebra $E=\roman{End}(T|_{D_F})$
with $D_F$ a finite diagram and $T$ a representation of a Nori diagram $D$,
the $R$--dual $E^\vee=\roman{Hom}_R(E,R)$ has the property that the canonical map
$E^\vee \otimes_R E^\vee \to \Hom(E,E^\vee)\simeq (E\otimes_R E)^\vee$
is an isomorphism. Thus, an $E$-module that is finitely generated projective as
an $R$-module carries the structure of an $E^\vee$-comodule.  The coalgebra
$\Cal{A}(D,T_\Z)$ also carries an algebra structure induced by the monoidal structure of
$\Cal{C}(D,T_\Z)$, see Sections 7.1.4 and 8.1 of [HuM-S17]. Thus, $\Cal{A}(D,T_\Z)$
determines a pro-algebraic monoid scheme $\roman{Spec}(\Cal{A}(D,T_\Z))$ (see 
Section 7.1.4 of [HuM-S17]).
This completes the proof.

\smallskip

In the more general setting of Section 4 below
we will only assume that the target category of our fiber functors is an abelian tensor
category, but not necessarily a category $R$-Proj with $R$ a Dedekind domain as
here above, with the Tannakian formalism of Sections 7.1.2-7.1.4 of [HuM-S17].
Indeed, the target category in general will be a category $\roman{Vec}^{(S,\leq)}$
for some poset $(S,\leq)$. This has a tensor structure obtained by 
identifying it with the category of covariant functors $\Cal{F}((S,\leq),\roman{Vec})$, 
endowed with the pointwise monoidal structure induced by the monoidal structure 
on $\roman{Vec}$. 

\bigskip

\centerline{2. THIN CATEGORIES AND PERSISTENCE}

\bigskip

In this section, we describe a slightly more general framework for persistence
constructions. It was sketched in the Introduction and its various
more precise versions will be given in the remaining Sections of the article.
\smallskip
We start with a category of geometric objects, and an indexing system for
persistence which is given by a thin category.

\medskip
{\bf 2.1. Geometric poset objects and thin categories.}
Let $\Cal{C}_{\roman{geom}}$ be a category of geometric objects (topological spaces, simplicial sets,
smooth manifolds, algebraic varieties, etc.). Whenever it is fixed, we refer to its morphisms
as ``geometric morphisms'' etc.

\smallskip

{\bf 2.1.1. Definition.} {\it
A poset object in $\Cal{C}_{\roman{geom}}$ is an object $S$ together with a subobject
$\Cal{R}\subset S\times S$ with the following properties:

$\bullet$ $(s,s)\in \Cal{R}$ for all $s\in S$,

$\bullet$  If $(s,s')\in \Cal{R}$ and $(s',s'')\in \Cal{R}$, then $(s,s'')\in \Cal{R}$;

$\bullet$ If $(s,s')\in \Cal{R}$ and $(s',s)\in \Cal{R}$, then $s=s'$.
\smallskip
The relation $(s,s')\in \Cal{R}$ we also denote by $s\leq s'$.
}

As was explained in Sec.~0.2 above, the notions of a poset  and of a thin category
essentially coincide.
\bigskip

{\bf 2.1.2. Remark.} The assumption that $\Cal{C}_{\roman{geom}}$ is a category
of geometric objects implies that points and subobjects are  defined in the usual
geometric terms. In a more general setting, one needs to use a formulation that
depends on a categorical notion of points in terms of the functor of points.
We will discuss this in the next section.

\medskip
{\bf 2.2. Persistence modules.} 
It is a general fact that the category of covariant functors $F:\Cal{B}\to \Cal{A}$
from a small category $\Cal{B}$ to an abelian category $\Cal{A}$ is itself abelian (see
for instance Proposition~44 of [Murf06]). Thus, we can give the following
definition.

\medskip
{\bf 2.2.1. Definition.} {\it Given a poset $(S,\leq)$ and an abelian category $\Cal{A}$, let $\Cal{A}^{(S,\leq)}$
be the abelian category whose objects  are the covariant functors $F: (S,\leq) \to \Cal{A}$ and morphisms
the natural transformations of such functors. Objects of $\Cal{A}^{(S,\leq)}$ will be referred to as
the $(S,\leq)$-persistence modules in $\Cal{A}$. In the case where $\Cal{A}=R$-Mod, we refer to
objects in $R$-Mod$^{(S,\leq)}$ as $(S,\leq)$-persistence $R$-modules. }

\bigskip

{\bf 2.3. Sublevel objects.}
As above let $\Cal{C}_{\roman{geom}}$ be a category of geometric objects and let
$(S,\leq)$ be a poset object in $\Cal{C}_{\roman{geom}}$ as in Definition 2.1.1.

\medskip

{\bf 2.3.1. Definition.} {\it
Consider a pair $(X, f:X \to S)$ where $X$ is an geometric object and $f$ its morphism
to the poset object $(S,\leq)$. For any $s\in S$ we define the sublevel objects $X_{f,s} \subset X$
as $X_{f,s}:=\{ x\in X\,:\, f(x)\leq s \in S \}$. They define inclusion maps
$j_{X,s,s'}: X_{f,s}\hookrightarrow X_{f,s'}$ for all $s\leq s'$ in $S$.}

\medskip

{\bf 2.4. Persistent functors.} Consider now a geometric object which is a poset
$(S,\leq)$ in its category, and let $\Cal{H}:\,\Cal{C}_{\roman{geom}}\to \Cal{A}$
be a functor with values in an abelian category $\Cal{A}$.

\smallskip

Denote by $\tilde\Cal{C}_{\roman{geom}}^S$  the 
category of ``geometric families over a  base $S$''. More precisely, one 
of its object
is a pair consisting of a geometric object and geometric morphism 
$(X,f:X\to S)$.  A morphism $(X,f) \to (X', f')$ is a geometric morphism 
$\varphi :\, X\to X'$ such that $f'\circ \varphi = f.$
 
 \smallskip
 
 Any functor $\Cal{H}: \Cal{C}_{\roman{geom}} \to \Cal{A}$ as above determines the functor
 $$
  \tilde\Cal{H}: \tilde\Cal{C}_{\roman{geom}}^S \to \Cal{A}^{(S,\leq)} 
 $$
 as follows. It sends each family $(X,f)$ to the object
$$
 \tilde\Cal{H}(X,f:X\to S) :=F_{X,f,S}: s \mapsto \Cal{H}(X_{f,s})
 $$
and each  morphism $\varphi : X\to X'$ of families  to the natural transformation of the functors 
$F_{X,f,S} \to F_{X',f',S}$ given by
$$
\tilde\Cal{H}(\varphi)  := \Cal{H}(\varphi_s): \Cal{H}(X_{f,s}) \to \Cal{H}(X'_{f',s}).
$$
Here $\varphi_s: X_{f,s} \to X'_{f',s}$ is the restriction $\varphi_s=\varphi|_{X_{f,s}}$
which maps to $X'_{f',s}$ because of the compatibility $f'\circ \varphi =f$.

\bigskip

{\bf 2.4.1. Definition.} {\it
For $\lambda\in S$, the persistence functor $P\Cal{H}_{\lambda}: \tilde\Cal{C}_{\roman{geom}}\to \Cal{A}^{(S,\leq)}$ of the
functor $\Cal{H}: \Cal{C}_{\roman{geom}} \to \Cal{A}$ is defined as follows.

It sends each family $f:X\to S)$ to
$$ 
P\Cal{H}_{\lambda}(X,f:X\to S):=s\mapsto \roman{Range}\,(\Cal{H}(j_{X_{s,\lambda}})) 
$$
for $s\leq \lambda$ and zero otherwise, where
$$
 j_{X_{s,\lambda}}: X_{f,s} \hookrightarrow X_{f,\lambda}, \ \ \ \text{ for } s\leq \lambda 
$$
and the induced morphism in $\Cal{A}$
$$ 
\Cal{H}(j_{X_{s,\lambda}}): \Cal{H}(X_{f,s}) \rightarrow \Cal{H}(X_{f,\lambda}). 
$$
On morphisms $\varphi: X\to X'$ with $f'\circ \varphi =f$ it is defined as the restriction of
$\tilde\Cal{H}(\varphi)$ to $\roman{Range}(\Cal{H}(j_{X_{s,\lambda}}))$ which takes values in
$\roman{Range}(\Cal{H}(j_{X'_{s,\lambda}}))$. }

\bigskip

{\bf 2.5. Example: persistent topology of graphs.}
Let $H$ be a finite directed graph of a thin category. Consider families of finite directed graphs
$(G,f:G\to H)$ over $H$. Let $G_{f,v}$ with $v\in V(H)$ be the respective  sublevel 
graphs. The graph $G_{f,v}$ is the induced
subgraph of $G$ on the set of vertices $w\in V(G)$ such that there is a path 
of directed edges in $H$ between $f(w)$ and $v$. The $(H,\leq)$-persistent 
topology of $G$ is then specified by the persistent connected components
$$
\roman{Range}\, ( H_0(G_{f,v};\Z) \to H_0(G_{f,v'};\Z)) \ \ \roman{for}\ v\leq v' \in V(H),
$$
and the persistent cycles
$$
\roman{Range}\, ( H_1(G_{f,v};\Z) \to H_1(G_{f,v'};\Z)) \ \ \roman{ for }\ v\leq v' \in V(H). 
 $$

 \medskip
{\bf 2.6. Example: persistent Orlik--Solomon algebras.}
Let $A$ be a hyperplane arrangement. Denote by $L(A)$ the associated
intersection poset, ordered by reverse inclusion. We consider morphisms $\varphi$
of hyperplane arrangements given by linear maps of the ambient space
that map one arrangement to the other and we write $L(\varphi)$ for the
induced map of intersection posets. We fix one arrangement $A$ and we consider
pairs $(B,\varphi)$ of arrangements endowed with a morphism $\varphi$ to $A$. 
The intersection poset $L(A)$ determines a structure of poset of topological
spaces (equivalently, an object in $\roman{Top}^{(L(A),\leq)}$ in the notation of [BuSc14],
[BuSiSc15]) on the hyperplane arrangement complement $M(A)$, with
inclusions $M(A)_s \hookrightarrow M(A)_{s'}$ for $s\leq s'$ in $L(A)$. Given a
morphism of arrangements $\varphi: B \to A$, defined as above, we obtain similarly a
structure of poset of topological spaces on the complement $M(B)$ indexed by
the poset $L(A)$. The families $M(B)_s$ with $s\in L(A)$ with the inclusions $M(B)_s
\hookrightarrow M(B)_{s'}$ for $s\leq s'$ in $L(A)$ form the system of sublevel objects
described above. The cohomology $H^*(M(B))$ (with coefficients in a field $\bold{K}$) 
of a hyperplane arrangement complement is the Orlik-Solomon algebra $OS(B)=H^*(M(B))$. 
We consider the cohomology $H^*(M(B)_s)$ 
with the maps induced by the inclusions $M(A)_s \hookrightarrow M(A)_{s'}$. In order to
have a covariant functor, we can consider the homology $OS^\vee(B):=H_*(M(B))$ with its structure
of module over the exterior algebra on $H^1(M(B))$, \cite{EPY03}. The associated
persistent functor is given by
$$ 
POS_\lambda^\vee(B):= \roman{Range}\,(H_*(M(B)_s) \to H_*(M(B)_\lambda)) , \ \ \roman{ for }\ s \leq \lambda \in L(A). 
$$

\medskip

Notice that we could also consider the persistent homology of a hyperplane arrangement
complement $M(B)$ with the persistence indexed by its own intersection poset $L(B)$.
The case considered above, where one considers arrangements $B$ endowed
with maps to a fixed arrangement $A$ and persistence with respect to the fixed $L(A)$
provides a uniform choice of the poset indexing the persistence modules. Allowing the
indexing poset to depend on the arrangement, as in the case where one uses $L(B)$
has advantages too, for example because there are in general few linear maps of
the ambient space that induce maps between two given arrangements.

\bigskip
{\bf 2.7. Example: persistent Tate motives.}
Beilinson, Goncharov, Schechtman and Varchenko in [BGSV90]  conjectured that, over 
a number field $\bold{K}$ the category of mixed Tate motive is generated by  motives of the form 
$$
\bold{m}\,(\bold{P}^n \smallsetminus A, A\smallsetminus (A\cap B)) 
 $$
where $A$ and $B$ are hyperplane arrangements in general position. These are
the motives whose cohomological realization gives the middle dimensional
relative cohomology of the pair $(\bold{P}^n \smallsetminus A, A\smallsetminus (A\cap B))$.
We can consider a setting as in the previous subsection, where one covers the hyperplane 
arrangement complement with sublevel sets indexed by a poset and correspondingly
consider motives of the form 
$\bold{m}\,((\bold{P}^n \smallsetminus A)_s, (A\smallsetminus (A\cap B)) _s)$ and persistent objects
$$
\roman{Range}\, (\bold{m}  ((\bold{P}^n \smallsetminus A)_s, (A\smallsetminus (A\cap B)) _s)\to \bold{m}\, 
((\bold{P}^n \smallsetminus A)_{s'}, (A\smallsetminus (A\cap B)) _{s'})) \   \roman{ for }\ s\leq s' .
 $$

\bigskip

\centerline{3. SUBLEVEL SIEVES AND PERSISTENCE}

\bigskip

In the setting described above we have assumed that we work with a category
$\Cal{C}_{\roman{geom}}$ of geometric objects and we have used the geometric
notion of points to define sublevel sets and persistence. We consider here
more general categories $\Cal{C}$ for which objects do not necessarily have
points in the geometric sense. However, they always have a ``functor of points" in Grothendieck's
sense: for an object $X$ in $\Cal{C}$ and another object $A$, an $A$-point
of $X$ is a morphism $\varphi: A \to X$ in $\roman{Mor}_\Cal{C}(A,X)$. We will use here
this approach to define a notion of persistent functors $P\Cal{H}$ associated to certain 
functors $\Cal{H}: \Cal{C} \to \Cal{A}$ with values in an abelian category. 

\bigskip

{\bf 3.1. Functor of points and Poset functors.}
Let $\Cal{C}$ be a category and $X\in \roman{Obj}\,(\Cal{C})$. The functor of points $\pi_X: \Cal{C} \to \roman{Sets}$
is a contravariant functor with $\pi_X(A) =\roman{Mor}_\Cal{C} (A,X)$ and $\pi_X(\varphi: B \to A)= - \circ \varphi:
\roman{Mor}_\Cal{C} (A,X) \to \roman{Mor}_\Cal{C}(B,X)$. 
The object $X$ is completely determined by its functor of points $\pi_X$
in the sense that a natural tranformation $\eta: \pi_X \to  \pi_Y$ determines a morphism $f: X \to Y$ in such a way 
that  natural equivalences are in (1,1)-correspondence with isomorphisms of the respective objects.

\medskip

{\bf 3.1.1. Definition.} {\it
Let $\Cal{C}$ be a category and $S$ its object. Let $\pi_S: \Cal{C} \to \roman{Sets}$ be the
functor of points of $S$. A poset functor on $S$ is
a contravariant functor $\Cal{R}_{(S,\leq)}: \Cal{C} \to \roman{Sets}$ given on objects by the assignment
of a subset $\Cal{R}_{(S,\leq)}(A)\subseteq \pi_S(A)\times \pi_S(A)$ with the following properties:

 (1) $(p_A,p_A)\in \Cal{R}_{(S,\leq)}(A)$ for all $p_A\in \pi_S(A)= \roman{Mor}_\Cal{C}(A,S)$;
 
(2) $(p_A,p_A')\in \Cal{R}_{(S,\leq)}(A)$ and $(p_A',p_A'')\in \Cal{R}_{(S,\leq)}(A)$ implies
$(p_A,p_A'')\in \Cal{R}_{(S,\leq)}(A)$;

(3)  $(p_A,p_A')\in \Cal{R}_{(S,\leq)}(A)$ and $(p_A',p_A)\in \Cal{R}_{(S,\leq)}(A)$ implies
$p_A=p'_A$ in $\pi_S(A)$;

(4)  if $(p_A,p_A')\in \Cal{R}_{(S,\leq)}(A)$ and $\varphi \in \roman{Mor}_\Cal{C}(B,A)$ then
$(p_A\circ \varphi, p_A'\circ \varphi)\in \Cal{R}_{(S,\leq)}(B)$.

The functor acts on morphisms by $\Cal{R}_{(S,\leq)}(\varphi: B \to A): \Cal{R}_{(S,\leq)}(A) \to \Cal{R}_{(S,\leq)}(B)$
mapping $(p_A,p_A')\mapsto (p_A\circ \varphi, p_A'\circ \varphi)$. }

\bigskip

{\bf 3.2. Sublevel sieve.} 
In this general setting, instead of the sublevel sets and sublevel objects we considered in the previous
sections, we construct sublevels as subfunctors of the functor of points.

\medskip

{\bf 3.2.1. Lemma.} {\it
(i) Let $\Cal{C}$ be a category with  terminal object $\star$. Let $(S,\leq)$ be a poset functor on
an object $S\in \roman{Obj}(\Cal{C})$ as in Definition 3.1.1.

Starting with a family $f:X\to S$ in $\Cal{C}$ and a 
choice of $s\in \roman{Mor}_\Cal{C}(\star,S)$, consider the assignment
$$ 
X_{f,\leq,s}(A) = \{ \alpha \in \pi_X(A)=\roman{Mor}_\Cal{C}(A,X)\,:\, f\circ\alpha \leq  s\circ t_A \} 
$$
where $t_A: A \to \star$ is the unique morphism to the terminal object in $\Cal{C}$ and 
$f\circ\alpha \leq  s\circ t_A$ means that $(f\circ\alpha, s\circ t_A)\in 
\Cal{R}_{(S,\leq)}(A) \subset \pi_S(A)\times \pi_S(A)$.

 Given a morphism
$\varphi: B \to A$ in $\roman{Mor}_\Cal{C}(B,A)$ assign to it the map
$$
 X_{f,\leq,s}(\varphi: B \to A) : X_{f,\leq,s}(A)  \to X_{f,\leq,s}(B),  \ \ \ \
 \alpha \mapsto \alpha \circ \varphi  . 
 $$
This assignment determines a contravariant functor $X_{f,\leq,s}: \Cal{C} \to \roman{Sets}$, 
the ``sublevel functor" of $f:X\to S$. 
\smallskip

(ii) For all $s\in \roman{Mor}_\Cal{C}(\star,S)$ the sublevel functor 
$X_{f,\leq,s}$ is a subfunctor of the functor of points $\pi_X$. Moreover, for 
$s \leq s'$  in $\roman{Mor}_\Cal{C}(\star,S)$, the sublevel functor 
$X_{f,\leq,s}$ is a subfunctor of $X_{f,\leq,s'}$. 

}

\medskip

{\bf Proof.}  (i) By Definition 3.1.1,  the subsets $\Cal{R}_{(S,\leq)}(A) \subset \pi_S(A)\times \pi_S(A)$
have the following property: if $(f\circ\alpha, s\circ t_A)\in \Cal{R}_{(S,\leq)}(A)$ then for any
$\varphi: B \to A$ in $\roman{Mor}_\Cal{C}(B,A)$ the element
$(f\circ\alpha\circ\varphi, s\circ t_A\circ\varphi)\in \Cal{R}_{(S,\leq)}(B)$, where
$t_A\circ\varphi=t_B$ is the unique morphism $t_B: B \to \star$ to the terminal object. Thus,
the assignment above is well defined and determines a contravariant functor. 
\smallskip
(ii) For all $A\in \roman{Obj}(\Cal{C})$ by construction we have $X_{f,\leq,s}(A)\subseteq \pi_X(A)$.
Moreover, for a morphism $\varphi: B\to A$ the image $X_{f,\leq,s}(\varphi)$ is the restriction
of $\pi_X(\varphi)$ (precomposition with $\varphi$) to $X_{f,\leq,s}(B)$. Hence $X_{f,\leq,s}$
is a subfunctor of the functor $\pi_X$. Similarly, if $s \leq s'$, that is $(s,s')\in \Cal{R}_{(S,\leq)}(\star)$,
the condition $f\circ \alpha \leq s \circ t_A$ implies that $f\circ \alpha \leq s' \circ t_A$ hence
$X_{f,\leq,s}(A)\subseteq X_{f,\leq,s'}(A)$ and $X_{f,\leq,s}(\varphi)$ is the restriction of
$X_{f,\leq,s'}(\varphi)$ hence $X_{f,\leq,s}$ is a subfunctor of $X_{f,\leq,s'}$. This completes the proof.

\smallskip
An assignment of a subfunctor of the functor of points $\pi_X$ is a sieve on $X$. Thus,
we equivalently refer to $X_{f,\leq,s}$ as the {\it sublevel sieve} of $X$.

\medskip

{\bf 3.2.2. Definition.} {\it
If the subfunctor $X_{f,\leq,s}$ of the functor of points $\pi_X$ is representable,
the object $X_s\in \roman{Obj}(\Cal{C})$ with $X_{f,\leq,s}(A)=\roman{Mor}_\Cal{C}(A,X_s)$
is the ``$s$-sublevel object" of $X$.}

\smallskip

Cases where the sublevel functor $X_{f,\leq,s}$ is representable include geometric
cases where it is a closed subfunctor. 

\smallskip

More precisely,  if the sublevel functor $X_{f,\leq,s}$ is representable by an object $X_s\in \roman{Obj}(\Cal{C})$,
then for any $s,s'\in \roman{Mor}_\Cal{C}(\star,S)$ with $s \leq s'$ there is a monomorphism
$j_{s,s'}: X_s \hookrightarrow X_{s'}$, since the inclusions $X_{f,\leq,s}(A)\subseteq X_{f,\leq,s'}(A)$
are monomorphisms of sets $j_{s,s'}: \pi_{X_s}\to \pi_{X_{s'}}$ which induce corresponding morphisms
$X_s \to X_{s'}$ in $\Cal{C}$ with the property that for all $u,v\in \pi_{X_s}$ if $j_{s,s'}\circ u=j_{s,s'}\circ v$
then $u=v$, hence $j_{s,s'}$ is a monomorphism in $\Cal{C}$. In the case of a representable sublevel
functor we can define  persistent functors in the following way.

\smallskip

{\bf 3.2.3. Definition.} {\it
Let $\Cal{C}$ be a category as above with  terminal object $\star$ and $(S,\leq)$
be a poset functor on an object $S\in \roman{Obj}(\Cal{C})$. Let $\Cal{H}: \Cal{C} \to \Cal{A}$ be
a covariant functor to an abelian category. For $f: X \to S$, and $s\in \roman{Mor}_\Cal{C}(\star,S)$
let $X_{f,\leq,s}$ be the sublevel functor. If $X_{f,\leq,s}$ is representable by $X_s\in \roman{Obj}(\Cal{C})$,
then the persistent functor $P\Cal{H}$ is given by
$$
P\Cal{H}_{(s,s')}(X) =\roman{Range}\,(\Cal{H}(j_{s,s'}): \Cal{H}(X_s) \to \Cal{H}(X_{s'})) \ \ \  \roman{ for }\ s\leq s' . 
$$
}

\smallskip

Return to the category $\tilde\Cal{C}_S$ of families in $\Cal{C}$.
For $(S,\leq)$ as above and $s\in \roman{Mor}_\Cal{C}(\star,S)$, and for an abelian category $\Cal{A}$ ,
we define $\Cal{A}^{(S,\leq)}$ as the category of covariant functors $F: \Cal{R}_{(S,\leq)}(\star) \to \Cal{A}$
and natural transformations of such functors. We can then interpret, for fixed $s'\in \roman{Mor}_\Cal{C}(\star,S)$
the persistent functor $P\Cal{H}$ as a functor $P\Cal{H}_{s'}: \tilde\Cal{C} \to \Cal{A}^{(S,\leq)}$ with
$$ 
P\Cal{H}_{s'}(X,f)= s\mapsto P\Cal{H}_{(s,s')}(X) 
$$
for $s\leq s'$ and zero otherwise. For a morphism $\varphi: (X,f) \to (X',f')$ 
we define $P\Cal{H}_{s'}(\varphi: (X,f) \to (X',f'))$ the map induced on
$\roman{Range}(\Cal{H}(j_{s,s'}): \Cal{H}(X_s) \to \Cal{H}(X_{s'}))$ by $\Cal{H}(\varphi): \Cal{H}(X_s) \to \Cal{H}(X'_s)$.

\bigskip

\centerline{4. NORI DIAGRAMS AND TANNAKIAN FORMALISM}

\medskip
{\bf 4.1. Setup.} Start with a diagram $D$ and an $R$-linear abelian category $\Cal{A}$, where $R$ is a commutative ring. 
Consider a representation
$T: D \to \Cal{A}^{(S,\leq)}$, where $(S,\leq)$ is a thin category, and $\Cal{A}^{(S,\leq)}=\Cal{F}((S,\leq),\Cal{A})$ 
is the category of covariant functors  from a thin category $(S,\leq)$ to 
$\Cal{A}$.
\smallskip
The representation $T$ assigns to each vertex
$v\in V(D)$ a functor $T(v): (S,\leq) \to \Cal{A}$ and to each edge $e\in E(D)$ a natural
transformation $T(e): T(s(e))\to T(t(e))$ between the functors associated to the source $s(e)$ and
target $t(e)$ vertices of $e$. 

\smallskip

For a vertex $v\in V(D)$, denote by  $\Cal{N}_v$ be the  set of natural self-transformations
$\alpha_v: T(v)\to T(v)$ of the functors $T(v): (S,\leq)\to \Cal{A}$.  Moreover,  put
$$
 \Cal{N}(T)=\{ (\alpha_v)_{v\in V(D)}\,:\, \alpha_v\in \Cal{N}_v \ \roman{ with }\
\alpha_{t(e)}\circ T(e)=T(e)\circ \alpha_{s(e)}, \, \forall e\in E(D)\}. 
$$

\smallskip

Generally, given an abelian category $\Cal{B}$ and a set $\Cal{S}$ of objects in $\Cal{B}$. As in
[HuM-S17],  denote by $\langle \Cal{S} \rangle$ the smallest full abelian
subcategory of $\Cal{B}$ that contains $\Cal{S}$ and such that the inclusion functor
is exact. It is generated by the objects in $\Cal{S}$ and is closed under
taking direct sums, direct summands, kernels and cokernels.
\medskip

Notice that here we do not assume that we start with 
a representation of the diagram $D$ in a category of $R$-modules for some
ring $R$. Hence we do not
have an obvious choice of a faithful exact functor from $\Cal{N}(T)$-Mod to
$\Cal{A}^{(S,\leq)}$ playing the role of the forgetful functor to $R$-Mod in the
setting of [HuM-S17]. However, we can still construct an abelian
category $\Cal{C}(D,T,\Cal{A}^{(S,\leq)})$ associated to the data $(D,T,\Cal{A}^{(S,\leq)})$ 
with the property that the representation $T: D \to \Cal{A}^{(S,\leq)}$ factors
through $\Cal{C}(D,T,\Cal{A}^{(S,\leq)})$.
\medskip

{\bf 4.1.1. Lemma.} {\it
Let  $T: D \to \Cal{A}^{(S,\leq)}$ be a diagram representation as above. Consider the
abelian subcategory $\langle T(D) \rangle$ of $\Cal{A}^{(S,\leq)}$. There is an
inclusion functor $\langle T(D) \rangle \hookrightarrow \Cal{N}(T)$-Mod. Let
$\Cal{C}(D,T,\Cal{A}^{(S,\leq)})$ denote the subcategory of $\Cal{N}(T)$-Mod
obtained in this way. If the inclusion functor $\langle T(D) \rangle \hookrightarrow \Cal{N}(T)$-Mod
is exact, this is an abelian subcategory. Moreover, there is a representation
$\tilde T: D \to \Cal{C}(D,T,\Cal{A}^{(S,\leq)})$ such that $T$ factors as $T=F\circ \tilde T$
with a faithful exact functor $F: \Cal{C}(D,T,\Cal{A}^{(S,\leq)}) \to \Cal{A}^{(S,\leq)}$.}

\medskip

{\bf Proof.} 
Since $\Cal{A}$ is an $R$-linear abelian category, the category $\Cal{A}^{(S,\leq)}$ 
is also an $R$-linear abelian category (Proposition~44 of [Mu06]), and the 
sets $\Cal{N}_v$ and $\Cal{N}(T)$ are rings with respect to the composition operation. 
We proceed as in Proposition~7.3.24 of [HuM-S17], 
assuming for simplicity that $D$ is a finite diagram. 

Consider the 
object $X=\oplus_v T(v)$. Then $\langle T(D) \rangle=\langle X \rangle$. Let
$\Cal{N}_X$ be the set of natural self-transformations of the functor $X: (S,\leq)\to \Cal{A}$. 
Among the transformations in $\Cal{N}_X$ we can identify
the elements of $\Cal{N}(T)$ as those $\alpha\in \Cal{N}_X$ that commute with the projections
$p_v: X\to T(v)$ and with the transformations $T(e)$. Thus, we can view
$\langle T(D) \rangle$ as a subcategory of $\Cal{N}(T)$-Mod. 
If the inclusion functor $\langle T(D) \rangle \hookrightarrow \Cal{N}(T)$-Mod
is an exact functor then $\langle T(D) \rangle$ gives an abelian subcategory of
$\Cal{N}(T)$-Mod, which we denote by $\Cal{C}(D,T,\Cal{A}^{(S,\leq)})$. There is 
a representation $\tilde T: D \to \Cal{C}(D,T,\Cal{A}^{(S,\leq)})$ given by the
assignment $v\mapsto T(v)$, $e\mapsto T(e)$ of the representation 
$T: D \to \Cal{A}^{(S,\leq)}$, seen as objects of $\langle T(D) \rangle$. By construction, this
representation satisfies  $T=F\circ \tilde T$, 
where $F: \Cal{C}(D,T,\Cal{A}^{(S,\leq)}) \to \langle T(D) \rangle \hookrightarrow \Cal{A}^{(S,\leq)}$
is the forgetful functor that forgets the $\Cal{N}(T)$-module structure.

This completes the proof.

\medskip

{\bf 4.1.2. Remark.} If the abelian category $\Cal{A}$ has a tensor structure and $\Cal{A}^{(S,\leq)}$
is endowed with the pointwise tensor structure, and we assume that the
diagram $D$ is a graded diagram with a commutative product with unit,
we can consider representations $T: D \to \Cal{A}^{(S,\leq)}$ that are unital
graded multiplicative representations so that $\Cal{C}(D,T,\Cal{A}^{(S,\leq)})$
also has a natural tensor structure such that $F: \Cal{C}(D,T,\Cal{A}^{(S,\leq)}) \to \Cal{A}^{(S,\leq)}$
becomes a tensor functor. (The argument for the original setting of 
diagram representations to $R$-Mod is given in Proposition~8.1.5(1) 
of [HuM-S17].) We can then consider faithful exact tensor functors 
from $\Cal{C}(D,T,\Cal{A}^{(S,\leq)})$ to categories $\Cal{B}^{(S,\leq)}$, where
$\Cal{B}$ is an $R$-linear abelian tensor category and $\Cal{B}^{(S,\leq)}$ is endowed with the pointwise
tensor structure. In particular, one can study tensor functors 
$\Cal{C}(D,T,\Cal{A}^{(S,\leq)}) \to \roman{Vec}^{(S,\leq)}$, as 
generalizations of persistent homology. 


\bigskip

\centerline{5. PERSISTENCE OF NORI MOTIVES}

\bigskip

We return here to the algebraic geometric environment  of Nori motives, in its 
updated form of Arapura's
category of Nori Motivic Sheaves, [Ar13].  We enrich it with a Persistence
structure. We write here the Nori motives covariantly (homologically) as in
[Ar13], rather than contravariantly (cohomologically) as in the initial
Grothendieck's project and in [HuM-S17].

\medskip

Let $S$ be a connected variety over $\bold{K}$.
The Nori Diagram $D_{ms}(S)$ for Motivic Sheaves over a base $S$ has the following structure (]Ar13]):
\smallskip
$\bullet$ One vertex in $V(D_{ms}(S))$ is a quadruple $(f: X \to S, Y, k, w)$, where
$X$ is a $\bold{K}$-variety with a morphism $f: X \to S$; $j: Y\hookrightarrow X$ is a closed
embedding (endowed with the restriction $f|_Y: Y \to S$), $k\in \Z_+$ is a non-negative
integer. and $w\in \Z$ is an integer.

\smallskip
$\bullet$ Edges in $E(D_{ms}(S))$ are of the following types:
\smallskip
(1) {\it Geometric morphisms}: each morphism of varieties $\varphi: X \to X'$ compatible with the maps to the base $S$ and the inclusions via commutative diagrams
$$ 
\xymatrix{ X \ar[rr]^\varphi \ar[dr]^f & & X' \ar[dl]_{f'} \\ & S & } \ \ \ \text{ and } \ \ \ 
\xymatrix{ Y \ar[r]^{\varphi |_Y} \ar[d]^j & Y' \ar[d]_{j'} \\ X \ar[r]^{\varphi} & X' }  
$$
produces an edge 
$$ 
(f: X\to S, Y, k, w) \to (f': X'\to S,Y',k,w) .
$$

\smallskip

(2) {\it Connecting morphisms}: every chain of closed embeddings $Z\subset Y \subset X$ determines
an edge
$$ 
(f: X \to S, Y, k+1,w) \to (f|_Y: Y \to S, Z, k, w) 
$$
(3) {\it Twist morphisms}: for every vertex $(f: X \to S, Y, k, w)$ there is an edge
$$ 
(f\circ p_1 : X\times \P^1\to S, Y\times \P^1 \cup X \times \{ 0 \}, k+2, w+1) \to (f: X \to S, Y, k, w) 
$$

A representation of the diagram $D_{ms}(S)$ is given by the constructible sheaves
$$
T_{ms}: (f: X \to S, Y, k, w) \mapsto H^k_S(X,Y,\F) . 
$$

The category of Nori Motivic Sheaves is defined as the Nori Diagram Category
$\Cal{C}(D_{ms}(S),T_{ms})$ of this representation. 

\bigskip

{\bf 5.1. Posets and semigroups.}
In order to enrich this construction with a
persistent version of the constructible sheaves $H^k_S(X,Y,\F)$, we follow the
general formalism described in the previous sections. The first step is to
consider a base $S$ that is endowed with a poset structure that will provide
the indexing of the persistence modules. 

\smallskip

A source of poset structures on geometric spaces such as $\bold{K}$-varieties is
the presence of a semigroup structure. Indeed there are natural poset structures
associated to semigroups, which we review briefly, see [Mi86]. 

\smallskip

Given a semigroup $S$, the set $E_S$ of idempotents $e\in S$, $e^2=e$, is
partially ordered by the relation $e\leq e'$ if $e=ee'=e'e$. A natural partial order
structure on a semigroup is one that restricts to this relation on the set of
idempotents. Recall that a semigroup $S$ is regular if every element $s\in S$
has a pseudoinverse $x$ such tat $sxs = s$.

\smallskip
The Nambooripad poset structure on a regular semigroup ([Na80]) 
is defined by 
$$ 
s \leq s' \ \ \ \roman{ iff } \ \ s = e s' = s' e' , \ \ \roman{ for\,\, some }\  e,e'\in E_S. 
$$
Several equivalent definitions are discussed in  [Mi86]. This is further
extended to more general semigroups as shown in [Mit86] as follows. 
Let $\hat S$ be obtained by adjoining a unit to the semigroup $S$ (or $S$ 
itself if it already has a unit). If $E_S$ is a subsemigroup of $S$ then the
relation $s\leq s'$ iff $s=es'=s'e'$ for some $e,e' \in E_{\hat S}$ is a partial
order compatible with multiplication. On an arbitrary semigroup the relation 
$$ 
s\leq s' \ \ \roman{ iff } \ \  s=xs'=s'y \ \ \roman{ and } \ \ 
xs=s \ \ \roman{ for\ some }\  x,y\in \hat S 
$$ is a partial order, which we call
as the natural partial order (see Proposition~2 and Theorem~3 of  [Mi86]). 

\bigskip

{\bf 5.2. Varieties with semigroup structures.}
In order to obtain a base $S$ with a poset structure $(S,\leq)$ it is then sufficient
to consider the case where $S$ has a semigroup structure. Semigroup structures
on algebraic varieties have been investigated in [Br14a]. We review here
some general facts from [Br14a] and [Br14b] and some 
examples from [Br14a].

\smallskip

We consider varieties over an algebraically closed field $\K$. Let $S$ be an
algebraic semigroup, that is, an algebraic variety over $\K$ endowed with
semigroup operation. Let $E_S$ be the subscheme of idempotents of $S$.
The partial order structure on idempotents in $E_S$ is given by $e_0\leq e_1$ iff 
$e_0 = e_0 e_1 = e_1 e_0$. Given two idempotents $e_0,e_1$ with $e_0\leq e_1$, 
the ``interval" $[e_0,e_1]$ is given by
$$
[e_0,e_1]:=\{ e\in E_S\,:\, e_0 \leq e \leq e_1 \} 
$$
and is a closed subscheme of $S$ ([Br14b], Corollary~2.17).
If $S$ is a smooth algebraic semigroup, then the scheme $E_S$ of idempotents
is also smooth ([Br14b], Remark~2.15). If $S$ is a commutative algebraic
semigroup, then the scheme $E_S$ of idempotents is finite and reduced 
([Br14b], Theorem~1.2). If $S$ is irreducible, then there is a smallest
closed and irreducible subsemigroup of $S$ that contains $E_S$, which is
given by a toric monoid ([Br14b], Theorem~1.2).
\smallskip
Classes of examples of algebraic varieties with semigroup structures include:

$\bullet$ linear algebraic semigroups: subsemigroups of $\roman{End}(V)$ with $V$
some finite dimensional vector space;
\smallskip
$\bullet$ an arbitrary variety $S$ with the left/right projection semigroup laws
$\mu_L(x,y)=x$ and $\mu_R(x,y)=y$ for all $x,y\in X$;
\smallskip
$\bullet$ algebraic groups;
\smallskip
$\bullet$  algebraic semigroup laws on Abelian varieties classified in Section 4 of [Br14a].

\smallskip
$\bullet$  algebraic semigroup laws on affine monomial curves (Theorem~5 of [Br14a]).

\bigskip

{\bf 5.3. Sublevel subschemes over semigroups.}
We consider the cases where $E_S$ is a subsemigroup of $S$ so that 
$S$ has a Nambooripad poset structure.

\smallskip

{\bf 5.3.1. Lemma.}{\it
Let $S$ be an algebraic semigroup, such that the idempotent subscheme $E_S$ is 
a subsemigroup of $S$. Let $(S,\leq)$ be the Nambooripad poset structure. 
For a morphism $f: X \to S$, the sublevel sets $X_s=\{ x\in X\,:\, f(x)\leq s \}$ are closed 
subschemes of $X$. }

\smallskip
{\bf Proof.} The subset $S_s=\{ a\in S\,:\, a \leq s \}$ is given by all the elements $a\in S$ that
are of the form $a=es =s e'$ for some pair $(e,e')\in E_S^2$. Thus, we have
$S_s = E_S \cdot s\, \cap\, s \cdot E_S \subset S$. For an algebraic semigroup, 
the subscheme $E_S$ of idempotents is closed, hence for a fixed $s\in S$ 
we obtain $S_a$ as an intersection of two closed subschemes.   
The sublevel sets $X_s=f^{-1}(S_s)$ are then preimages in $X$ of 
closed subschemes of $S$.
This completes the proof.

\smallskip

The sublevel subvarieties $X_s$ have embeddings $j_{s,s'}: X_s\hookrightarrow X_{s'}$ for
$s\leq s'$ due to the transitive property of the partial order relation on $S$. 

\smallskip

{\bf 5.3.2. Corollary.}{\it
If $S$ is a commutative algebraic semigroup, which is smooth as an algebraic variety, 
and $(S,\leq)$ is the Nambooripad poset structure, then the sublevels $X_s$ 
of a morphism $f: X \to S$ are finite unions of fibers. }
\smallskip
{\bf Proof.} For a commutative smooth algebraic semigroup $S$, the idempotent subscheme $E_S$
is a subvariety of $S$ consisting of a finite set of points. Then the $S_s$ are also finite and identified
with $\{ es\, e\in E_S \}$, and the preimages $X_s=f^{-1}(S_s)=\cup_{e\in E_S} f^{-1}(es)$ are finite
unions of fibers of the morphism $f: X \to S$.

\bigskip

{\bf 5.4. Persistent Nori Motivic Sheaves.}
Let $S$ be a $\K$-variety with a semigroup structure and an associated
partial order $(S,\leq)$ as discussed above. As in Arapura's setting of Nori Motivic Sheaves 
in [Ar13] we consider the Nori Diagram with vertices $(f:X\to S, Y, k,w)$
and with the three classes of edges described above. To this diagram $D_{ms}(S,\leq)$
we associate a representation in the abelian category $\roman{Vec}^{(S,\leq)}$ in the 
following way. 

\smallskip

We consider a Persistence Diagram $D^P_{ms}(S,\leq)$ of Motivic Sheaves over a base
algebraic semigroup $S$ with an associated partial order $\leq$ determined by the
semigroup structure with the following vertices and edges:
\smallskip

$\bullet$  Vertices in $V(D^P_{ms}(S,\leq))$ are given by elements of the form $(f: X \to S, Y, k, w, s)$, 
 where $X$ is a $\K$-variety with a morphism $f: X \to S$, $j: Y\hookrightarrow X$ is a closed
embedding (endowed with the restriction $f|_Y: Y \to S$), $k\in \Z_+$ is a non-negative
integer, $w\in \Z$ is an integer, and $s\in S$.

\smallskip
$\bullet$ Edges in $E(D^P_{ms}(S,\leq))$ are of the following types:
\smallskip
(1) {\it Geometric morphisms}: for a morphism of varieties $\varphi: X \to X'$ compatible with the maps to the base $S$ and the inclusions via commutative diagrams
$$
 \xymatrix{ X \ar[rr]^\varphi \ar[dr]^f & & X' \ar[dl]_{f'} \\ & S & } \ \ \ \text{ and } \ \ \ 
\xymatrix{ Y \ar[r]^{\varphi |_Y} \ar[d]^j & Y' \ar[d]_{j'} \\ X \ar[r]^{\varphi} & X' }  
$$
there corresponds an edge 
$$ 
(f: X\to S, Y, k, w, s) \to (f': X'\to S,Y',k,w, s)
 $$
 \smallskip
(2) {\it Connecting morphisms}: every chain of closed embeddings $Z\subset Y \subset X$ determines
an edge
$$ 
(f: X \to S, Y, k+1,w, s) \to (f|_Y: Y \to S, Z, k, w, s) 
$$
\smallskip

(3) {\it Twist morphisms}: for every vertex $(f: X \to S, Y, k, w)$ there is an edge
$$
 (f\circ p_1 : X\times \P^1\to S, Y\times \P^1 \cup X \times \{ 0 \}, k+2, w+1, s) \to 
(f: X \to S, Y, k, w, s) 
$$

\smallskip

(4) {\it Persistence morphisms}: for all $s'\in S$ with $s \leq s'$  an edge
$$
(f: X\to S, Y, k, w, s) \to (f: X\to S, Y, k, w, s'). 
$$
\medskip
A representation of the diagram $D^P_{ms}(S,\leq)$ in $\roman{Vec}^{(S,\leq)}$
is then obtained as follows. Consider the map that assigns to each vertex
$(f: X \to S, Y, k, w,(s,s'))$ the functor
$$ F: (S,\leq) \to \roman{Vec}, $$
where we view $(S,\leq)$ as a thin category, given by
$$ 
F(t) = \roman{Range}(H_*(j_{t,s}): H_*(X_t,Y_t;\Q) \to H_*(X_s,Y_s;\Q)), 
$$
for all $t\leq s$ in $S$ and zero otherwise. The edges listed above are correspondingly
mapped to homomorphisms of the homology groups $H_*(X_t,Y_t;\Q)$ and $H_*(X_s,Y_s;\Q)$
with induced morphisms on $\roman{Range}(H_*(j_{t,s}))$.

\bigskip

\centerline{6. MODEL CATEGORIES AND PERSISTENT TOPOLOGY}
\bigskip
In this section we address a question that was posed to us by Jack Morava,
about developing a suitable model structure for persistent topology. 

\medskip
{\bf 6.1. Model categories.}
We review quickly some basic definitions regarding model structures and
categories that we will be using in the following. 

\smallskip

Recall that a morphism $f$ in a category is called a retract of another morphism $g$ iff there is a
commutative diagram
$$ 
\xymatrix{ A \ar[d]^f \ar[r] & C \ar[r] \ar[d]^g & A \ar[d]^f \\ B \ar[r] & D \ar[r] & B } 
$$
where the horizontal compositions are the identities.

\smallskip

A model category $\Cal{M}$ is a category together with three classes of morphisms: {\it weak equivalences, fibrations}
 and {\it cofibrations.} These data must satisfy the following axioms:

\medskip

(1) $\Cal{M}$ has all small limits and colimits;
\smallskip
(2) if in a composition $g\circ f$ of morphisms two among the three maps $f,g,g\circ f$
are weak equivalences, then the third also is;
\smallskip
(3) if a map $f$ is a retract of $g$, and if $g$ is a weak equivalence, fibration, or cofibration,
then $f$ also is;
\smallskip
(4) given a commutative diagram
$$ 
\xymatrix{ A \ar[r] \ar[d]^\iota & X \ar[d]_p \\ B \ar[r]&  Y } 
$$
a lift $B\to X$ exists if either $\iota$ is a cofibration and $p$ is an acyclic fibration (i.e.~both a fibration
and a weak equivalence) or if $\iota$ is an acyclic cofibration (i.e.~both a cofibration
and a weak equivalence) and $p$ is a fibration;
\smallskip
(5) morphisms $g$ in $\Cal{M}$ can be factored as $g=q i$ with $q$ an acyclic fibration and $i$ a cofibration or as $g=pj$ with $p$ a fibration and $j$ an acyclic cofibration.

\medskip

Let $\Cal{M}, \Cal{N}$ be model categories. {\it A Quillen pair $L: \Cal{M} \leftrightarrow \Cal{N}: R$} is an adjoint
pair of functors $(L,R)$ where $L$ preserves cofibrations and $R$ preserves fibrations, see Section~1.6
of [Be18].

\smallskip

A model category $\Cal{M}$ is {\it cofibrantly generated} iff:
\smallskip
(1) there is a set $\Cal{}I$ of cofibrations of $\Cal{M}$
such that the domains of the elements of $\Cal{I}$ are small with respect to $\Cal{I}$  and such that a map is an acyclic fibration iff it has the right lifting property
with respect to $\Cal{I}$;
\smallskip
(2) there is a set $\Cal{J}$ of acyclic cofibrations of $\Cal{M}$ (also with the property
that the domains are small with respect to $\Cal{J}$) such that a map is a fibration iff it has the right lifting property with respect to 
$\Cal{J}$.

\medskip
The set $\Cal{I}$ is the set of generating cofibrations and the set $\Cal{J}$ is the set of generating
acyclic cofibrations
(see Section 1.7
of [Be18] for the terminology and for more details).

\smallskip

A model category $\Cal{M}$ is called {\it combinatorial} if it is cofibrantly generated and as a 
category it is locally presentable, that is, it admits all small colimits and a set of
small objects such that any object can be obtained as colimit of a small diagram 
with objects in this set, see Sec. 2.7 of [Be18] for a more detailed exposition.

\bigskip
{\bf 6.2. Model structures on categories of functors.}
Given a small category $\Cal{C}$ and a category $\Cal{D}$, we denote by $\Cal{F}(\Cal{C},\Cal{D})$ 
the category of functors, whose objects are the covariant functors $F:\Cal{C}\to \Cal{D}$
and morphisms the natural transformations of these functors. 

\smallskip

For example, the categories we considered in the previous sections of the form
$\roman{Vec}^{(S,\leq)}$ with $(S,\leq)$ a thin category
are categories of functors, and so are the $\roman{Top}^{(S,\leq)}$ considered
in [BuSc14], [BuSiSc15] as a setting for persistent topology.

\smallskip

It is known that if the category $\Cal{D}$ has a model structure that is cofibrantly
generated, then for any small category $\Cal{C}$ the category of functors $\Cal{F}(\Cal{C},\Cal{D})$
also has a model structure, called the {\it projective model structure}, which is also 
cofibrantly generated, see Sec. 11.6 of [HI03].
If the category $\Cal{D}$ has a combinatorial model structure, then the category of functors
 $\Cal{F}(\Cal{C},\Cal{D})$ has a model structure, called the {\it injective model structure}.

\smallskip

In the projective model structure on  $\Cal{F}(\Cal{C},\Cal{D})$
 weak equivalences and projective fibrations are those natural transformations 
$\eta: F\to F'$ of functors $F,F': \Cal{C} \to \Cal{D}$ such that for all objects $X\in \roman{Obj}(\Cal{C})$,
the morphisms $\eta_X: F(X)\to F'(X)$ in $\Cal{D}$ are, respectively, weak equivalences
and fibrations.  Similarly, in the injective model structure injective
weak equivalences and injective cofibrations are natural transformations that are,
object-wise in $\Cal{C}$, weak equivalences and cofibrations in~$\Cal{D}$.

\bigskip
{\bf 6.3. Model structure for datasets and Vietoris--Rips complexes.}
We consider here the main example of persistent topology, which accounts for its use
in topological data analysis, namely datasets with their associated Vietoris--Rips complexes,
and their persistent homology barcode diagrams. 

\smallskip

We construct a model category for this setting in several steps:
\medskip
(1) We start by considering the model structure on the categories of simplicial sets $\Delta\Cal{S}$ and
of chain complexes $\roman{Ch}_R$  over a commutative ring $R$. We denote by $\Cal{M}$ either of
these model categories.
\smallskip

(2) We induce a projective model structure on the category of functors 
$$
\Cal{M}^{(S,\leq)}=\Cal{F}((S,\leq),\Cal{M}),
$$
where $(S,\leq)$ is a suitable poset, viewed as a thin category.
\smallskip

(3) We construct a category $\Cal{P}_\bold{E}$ of finitely supported probability distributions in a fixed 
ambient metric space $\bold{E}$ (e.~g.~an Euclidean space of sufficiently large dimension).
\smallskip

(4) We describe the assignment of Vietoris--Rips complexes to datasets as a functor 
$VR: \Cal{P}_\bold{E} \to \Cal{M}^{(S,\leq)}$.
\smallskip

(5) We use Dugger's construction [Du01] of a universal model category associated to a
small category to construct a model category $U(\Cal{P}_\bold{E})$ for finitely 
supported probability distributions with a Vietoris--Rips functor 
$VR: U(\Cal{P}_\bold{E})\to \Cal{M}^{(S,\leq)}$.

\medskip
{\bf 6.3.1. Model structure on chain complexes.}
Let $\roman{Ch}_R$ be the category of (unbounded) chain complexes over a commutative ring $R$
with chain maps as morphisms. The weak equivalences are the quasi-isomorphisms of chain
complexes. The fibrations are the chain maps $\varphi_\bullet : C_\bullet \to C'_\bullet$
such that for each level $n$ the map $\varphi_n: C_n \to C'_n$ is an epimorphism of
$R$-modules.  The cofibrations are chain maps that are level-wise monomorphisms of $R$-modules
with projective cokernel.
\smallskip
This is a projective model structure on $\roman{Ch}_R$. One can similarly 
consider the injective model structure with the same weak equivalences, but with cofibrations
given by chain maps that are level-wise injective morphisms of $R$-modules, while
fibrations are level-wise epimorphisms with injective kernel. 
\smallskip
These two model structures are Quillen-equivalent (Section~1.7 of [Be18]).
The projective model structure on $\roman{Ch}_R$ is cofibrantly generated
(see Sec.~ 2.3.11 of [Hov98] and Sec. 5 of [SchSh00]).

\medskip
{\bf 6.3.2. Model structure on simplicial sets.}
The category of simplicial sets $\Delta\Cal{S}$ has a model structure (the Kan-Quillen 
model structure) where the weak equivalences are morphisms that induce a weak homotopy
equivalence of topological spaces at the level of geometric realizations, the fibrations are Kan
fibrations, and the cofibrations are monomorphisms of simplicial sets. The Kan--Quillen 
model structure is cofibrantly generated, with generating cofibrations the boundary inclusions
and generating acyclic cofibrations the horn inclusions, see [Be18], [GeMa03], 
[GoJa99], [Hi03]. 

\medskip
{\bf 6.3.3. Model structure on indexed diagrams.}
Let $\Cal{M}$ denote either the category of chain complexes $\roman{Ch}_R$ with the
projective model structure, or the category of simplicial sets $\Delta\Cal{S}$ with the 
Kan--Quillen model structure. We consider now the category $\Cal{M}^{(S,\leq)}$
of $(S,\leq)$-indexed diagrams in $\Cal{M}$, for a poset $(S,\leq)$. This is the same as the category 
$$ \Cal{M} =\Cal{F}((S,\leq),\Cal{M}) $$
of covariant functors from the thin category $(S,\leq)$ to $\Cal{M}$ with morphisms
given by natural transformations. This will include the case of $\Cal{M}^{(\R,\leq)}$
that will correspond to the usual Vietoris--Rips complexes of data sets, but we will
work with a more general $(S,\leq)$ that also incorporates in the Vietoris--Rips complex
a cutoff according to a probability (see Part~IV of  [BoChYv18]).

Since the model category $\Cal{M}$ is cofibrantly generated, the category of functors
$\Cal{M}^{(S,\leq)}$ admits a projective model structure that is also cofibrantly generated. 

\medskip

{\bf 6.3.4. A category of data sets.} 
Consider the following small category $\Cal{P}_\bold{E}$. Its objects are triples $(A,f,P)$ of
a finite set $A$, an embedding $f: A \hookrightarrow \bold{E}$ in a fixed ambient
metric space $\bold{E}$, which we can assume to be an Euclidean space of some fixed
sufficiently large dimension, and a probability distribution $P$ on $A$, which
we can view as a probability on $\bold{E}$ supported on the finite set $f(A)$ through
pushforward along the map $f$. 
\smallskip
A morphism in $\roman{Mor}_{\Cal{P}_\bold{E}}((A,f,P),(A',f',P'))$ 
is a 
pair $(\varphi, \tilde\varphi)$ consisting of 
a continuous and Lipschitz map $\tilde\varphi : \bold{E} \to \bold{E}$ 
that restricts to a map $\varphi: A \to A'$ through a commutative diagram
$$ 
\xymatrix{ A \ar[r]^{\varphi} \ar[d]^f & A' \ar[d]_{f'} \\ \bold{E} \ar[r]^{\tilde\varphi} & \bold{E} } 
$$
such that $P'=\varphi_* P$, the pushforward measure given by 
$$ 
(\varphi_* P)_y = \sum_{x\in \varphi^{-1}(y)} P_x,  \ \ \ \forall y \in A'. 
$$
Thus morphisms in $\Cal{P}_\bold{E}$ are subsets of the set of Lipschitz functions of $\bold{E}$.
The probability distribution $P$
on the finite set $A$ should be thought of as assigning a degree of reliability to the points
in the dataset, see the discussion in [BoChYv18]. Points $x\in A$ with a low
probability $P_x$ should be regarded as errors in the data and discarded in the 
construction of the associated simplicial complexes.

\bigskip
{\bf 6.3.5. Vietoris--Rips functors.}
Given an object $(A,f,P)$ in $\Cal{P}_\bold{E}$, we construct a Vietoris--Rips complex $VR(A,f,P)$,
obtained by considering, for any choice of an error threshold 
$\Lambda \in [0,1]$, the set of points
$$ 
X_{\Lambda}=\{ x\in f(A)\subset \bold{E}\,|\, P_x\geq \Lambda \} 
$$
and then constructing the Vietoris--Rips complex $VR_\bullet(X_\Lambda,t)$
for $t\in \R_+^*$ where $VR_n(X_\Lambda,t)$ is the span of all unordered
$(n+1)$-tuples of points $(x_0,\ldots, x_n)$ in $X_\Lambda$ for which all
the pairwise distances satisfy $\roman{dist}(x_i,x_j)\leq t$. The threshold $\Lambda$
is aimed at discarding a number of outliers of small probability among the data.

\smallskip

Consider then a morphism $(\varphi, \tilde\varphi)$ in $\Cal{P}_\bold{E}$, where $\tilde\varphi$
is a Lipschitz function $\tilde\varphi: \bold{E} \to \bold{E}$ with Lipschitz constant $K>0$ 
and $\varphi: A \to A'$ has $m=\min_{y\in A'} \# \varphi^{-1}(y)$. 
This map sends $X_\Lambda$ to $X'_{m\Lambda}=\{ y\in f'(A')\,|\, P_y\geq m\Lambda \}$.
Moreover, it sends a pair of points $x_i,x_j\in X_\Lambda$ with distance 
$\roman{dist}(x_i,x_j)\leq t$ to a pair of points $\tilde\varphi(x_i), \tilde\varphi(x_j)\in X'_{m\Lambda}$
with distance $\roman{dist}(\tilde\varphi(x_i),\tilde\varphi(x_j))\leq Kt$. Thus, it induces a morphism
$$ VR(\varphi,\tilde\varphi): VR_\bullet(X_\Lambda,t) \to VR_\bullet(X_{m\Lambda},Kt).$$

\smallskip

Consider the set $S=\R \times [0,1]$ with the partial order structure $\leq$ given by the
product order $(t,\Lambda)\leq (t',\Lambda')$ iff $t\leq t'$ and $\Lambda \geq \Lambda'$
in the natural ordering of the real numbers (with the reverse ordering on $[0,1]$). 
We regard $(S,\leq)$ as a thin category.
\medskip

{\bf 6.3.6. Proposition.} {\it
The assignments $(A,f,P)\mapsto VR(A,f,P)$ and $(\varphi,\tilde\varphi)\mapsto 
VR(\varphi,\tilde\varphi)$ as above determine a functor $\Cal{P}_\bold{E} \to \Cal{M}^{(S,\leq)}$.}

\smallskip
{\bf Proof.}
The inclusions $X_\Lambda \subset X_\Lambda'$ for $\Lambda \geq \Lambda'$
and the inclusions of the subset of points of $X_\Lambda$ with mutual distances
bounded above by $t$ in the subset of points with mutual distances at most $t'$
for $t\leq t'$ induce morphisms 
$$ j_{(t,\Lambda),(t',\Lambda')}: VR_\bullet(X_\Lambda,t) \to VR_\bullet(X_{\Lambda'},t') $$ 
for $(t,\Lambda)\leq (t',\Lambda')$ in the chosen ordering of $S$. Thus, the assignment
$$ VR(A,f,P): (t,\Lambda) \mapsto VR^{(t,\Lambda)}(A,f,P) := VR_\bullet(X_\Lambda,t)  $$
determines a functor $(S,\leq)\to \Cal{M}$, that is, an object in $\Cal{M}^{(S,\leq)}$.
The morphisms $VR(\varphi,\tilde\varphi)$ of Vietoris--Rips complexes 
described above, for 
$$
(\varphi,\tilde\varphi)\in \roman{Mor}_{\Cal{P}_\bold{E}((A,f,P),(A',f',P'))},
$$ 
determine natural transformations $\eta_{(\varphi,\tilde\varphi)}:
VR(A,f,P)\to VR(A',f',P')$ of the functors $VR(A,f,P), VR(A',f',P'): (S,\leq) \to \Cal{M}$,
with 
$$ \eta_{(\varphi,\tilde\varphi), (t,\Lambda)}=VR(\varphi,\tilde\varphi) : 
VR_\bullet(X_\Lambda,t)\to VR_\bullet(X_{m\Lambda},Kt) $$
satisfying 
$$ \eta_{(\varphi,\tilde\varphi), (t',\Lambda')}\circ j_{(t,\Lambda),(t',\Lambda')} = j_{(Kt,m\Lambda),(Kt',m\Lambda')} \circ \eta_{(\varphi,\tilde\varphi), (t,\Lambda)}, $$
for all $(t,\Lambda)\leq (t',\Lambda')$ in $(S,\leq)$.
\smallskip
This completes the proof.

\smallskip

In the special case where we fix $\Lambda=0$, hence we consider all the points of $f(A)\subset \bold{E}$
regardless of the assigned probabilities, then this construction recovers the usual Vietoris-Rips
complexes as objects in $\Cal{M}^{(\R,\leq)}$ and the functor $VR: \Cal{P}_\bold{E} \to \Cal{M}^{(\R,\leq)}$
factors through the forgetful functor $\Cal{P}_\bold{E} \to \Cal{D}_\bold{E}$ where $\Cal{D}_\bold{E}$ is the
category of unweighted data sets with objects $(A,f)$ and morphisms given by restrictions
of Lifschitz functions of $\bold{E}$ (with no conditions on probabilities). The Vietoris-Rips
functor then defines a functor $VR: \Cal{D}_\bold{E} \to \Cal{M}^{(\R,\leq)}$. It is convenient to include
the probability data in the construction. We will discuss a more general way of including
probability data in the next section. 

\medskip

{\bf 6.3.7. Dugger's universal model structure.}
Dugger's construction in [Du01] assigns a universal model category $U(\Cal{C})$ 
to a small category $\Cal{C}$, with the property that functors from $\Cal{C}$ to a model 
category factor through $U(\Cal{C})$. The main idea is that $U(\Cal{C})$ extends the
category $\Cal{C}$ by formally adjoining homotopy colimits. A factorization of
a functor $F: \Cal{C} \to \Cal{M}$, where $\Cal{M}$ is a model category, through another
model category $\tilde\Cal{M}$ with a functor $J: \Cal{C} \to \tilde\Cal{M}$ consists of a
Quillen pair $L: \tilde\Cal{M} \rightleftarrows \Cal{M} : R$ and a natural weak equivalence
$\eta: L\circ J \to F$. (For a brief review of Quillen pairs see Definition~1.6.3 and
Proposition~1.6.4 of [Be18].) The main result of [Du01]  shows that
given a small category $\Cal{C}$, there exists a closed model category $U(\Cal{C})$
with a functor $J: \Cal{C} \to U(\Cal{C})$ such that any functor $F: \Cal{C} \to \Cal{M}$ to a 
model category $\Cal{M}$ factors, in the sense recalled above, through $U(\Cal{C})$.

\medskip
{\bf 6.3.8. Proposition.} {\it
Let $VR: \Cal{P}_\bold{E} \to \Cal{M}^{(S,\leq)}$ be the Vietoris-Rips functor of Proposition 6.3.6,
where $\Cal{M}^{(S,\leq)}$ has the projective model structure. There is a Quillen pair
$L: U(\Cal{P}_\bold{E}) \rightleftarrows \Cal{M} : R$ and a natural weak equivalence
$\eta: L\circ J \to VR$ that factor the Vietoris--Rips functor
through the universal model category $U(\Cal{P}_\bold{E})$ of the category $\Cal{P}_\bold{E}$ of data sets.} 

\smallskip
{\bf Proof.}
We apply the construction of [Du01] to the category $\Cal{P}_\bold{E}$ and we obtain
an associated universal model category $U(\Cal{P}_\bold{E})$. The factorization property
discussed above implies that the Vietoris-Rips functor factors through $U(\Cal{P}_\bold{E})$,
when we consider $\Cal{M}^{(S,\leq)}$ endowed with a model structure. We have seen
above that $\Cal{M}^{(S,\leq)}$ always supports the projective model structure, seen
as the category of functors $\Cal{F}((S,\leq),\Cal{M})$.

\bigskip
\centerline{7. PERSISTENCE AND $\Gamma$-SPACES}
\medskip
{\bf 7.1. Setup.} We discuss in this section how to adapt to the context of persistent topology another
important homotopy-theoretic construction: Segal's $\Gamma$-spaces. 
We consider a (slightly modified) setting developed in [Mar18] that incorporates probabilistic data
in the construction of Segal's $\Gamma$-spaces. They replace  the
finite probability distributions considered in the previous section in the context of
persistent topology of databases.

\smallskip

Our point of view here  is somewhat different from the one proposed in [Mar18] and
more closely related to our discussion of the Vietoris-Rips functor in the previous section. The
point we want to stress here is that the Vietoris-Rips functor, as we described it above, can be
generalized using Segal's $\Gamma$-spaces ([Se74]). A $\Gamma$-space is a functor
$F: \Gamma^0 \to \Delta\Cal{S}_*$
from the category of pointed finite sets to pointed simplicial sets. In particular, Segal showed in
[Se74] that to a category $\Cal{C}$ with a categorical sum and a zero object one can associate
a $\Gamma$-space obtained by assigning to a finite pointed set $(X,\star)$ the nerve
$\Cal{N}\Sigma_\Cal{C}(X,\star)$ of the category $\Sigma_\Cal{C}(X,\star)$ of summing functors 
$\Phi: P(X,\star)\to\Cal{C}$, where $P(X,\star)$ is the category
with objects the pointed subsets of $X$ and morphisms the pointed inclusions, and the summing
functors satisfy $\Phi_{X,\star}(S)\oplus \Phi_{X,\star}(S') \simeq \Phi_{X,\star}(S\cup S')$ for all
$S,S'\in P(X,\star)$ with $S\cap S'=\{ \star \}$. A map of pointed sets $f: (X,\star)\to (Y,\star')$ induces on 
summing functors a transformation $\Sigma_\Cal{C}(f): \Sigma_{\Cal{C}}(X,\star)\to \Sigma_{\Cal{C}}(Y,\star')$ given
by $\Sigma_\Cal{C}(f)(\Phi_{X,\star})(S)=\Phi_{X,\star}(f^{-1}(S))$, for all $S\in P(Y,\star')$.

\smallskip 

As in the previous section, we consider a category of databases, identified with 
finite sets endowed with probabilities, embedded in an ambient metric space (a large
dimensional Euclidean space $\bold{E}$). To adapt the setting to the pointed case required
for the $\Gamma$-space formalism, we consider here pointed sets, so we work with
the category $\Cal{P}_{\bold{E},*}$ whose objects are triples $((X,\star),f,P)$ of a finite pointed
set $(X,\star)$ with a probability measure $P$ and an embedding $f: X \hookrightarrow \bold{E}$
with $f(\star)=0$ the origin, in the Euclidean space $\bold{E}$. Morphisms are pairs $(\varphi,\tilde\varphi)$
of a measure-preserving pointed map $\varphi$ of finite sets and a 
Lipshitz self-map $\tilde\varphi$ of $\bold{E}$ that fixes the origin, which restricts to $\varphi$
on the images under the embeddings.

\smallskip

Following the construction of the Vietoris-Rips complex in the previous section, 
we consider the poset $(S,\leq)$ with $S=\R_+\times [0,1]$ with the natural 
order on $t\in\R_+$ and the reverse order on $\Lambda\in [0,1]$.  We then
have the following generalizations of the Vietoris-Rips functor with values
in $\Delta\Cal{S}^{(S,\leq)}$.

\medskip

{\bf 7.2. Proposition.} {\it
Any $\Gamma$-space $F_\Cal{C}: \Gamma^0 \to \Delta\Cal{S}_*$ determines a functor}
$$
\tilde F_\Cal{C}: \Cal{P}_{\bold{E},*} \to \Delta\Cal{S}_*^{(S,\leq)}. 
$$
\smallskip
{\bf Proof.}
Start with a $\Gamma$-space $F_\Cal{C}: \Gamma^0 \to \Delta\Cal{S}_*$ associated as above to
a category $\Cal{C}$ with sum and zero object. 
\smallskip
We obtain from it a functor $\tilde F_\Cal{C}: \Cal{P}_{\Cal{E},*} \to
\Delta\Cal{S}_*^{(S,\leq)}$ in the following way. For $s=(t,\Lambda) \in S$ we denote by $(X_{t,\Lambda},\star)$
the finite pointed set $\{ \star \} \cup Y_{t,\Lambda}$ where $Y=X\smallsetminus \{ \star \}$ and
$Y_{t,\Lambda}$ is the subset of all the non-marked points of $X$ with mutual distances
$d_{\bold{E}}(f(x),f(x'))\leq t$ and with probabilities $P_x\geq \Lambda$, as in the previous section.
We can then set $\tilde F_\Cal{C} ((X,\star),f,P)(s) = F_\Cal{C}(X_s)$. For $s\leq s'$ in $S$ we have
an inclusion $X_s\hookrightarrow X_{s'}$. This induces a transformation as above $\Sigma_\Cal{C}(X_s)\to
\Sigma_\Cal{C}(X_{s'})$ on summing functors, hence a map $F_\Cal{C}(X_s)\to F_\Cal{C}(X_{s'})$, hence
we obtain a functor $\tilde F_\Cal{C}((X,\star),f,P): (S,\leq) \to \Delta\Cal{S}_*$. Morphisms
$(\varphi,\tilde\varphi): ((X,\star),f,P)\to ((X',\star'),f',P')$ in $\Cal{P}_{\bold{E},*}$ induce by restriction 
maps $\varphi_s: X_s \to X_{s'}$ and corresponding maps on summing functors $\Sigma_\Cal{C}(X_s)\to
\Sigma_\Cal{C}(X_{s'})$ as above.

\bigskip
\centerline{8. FURTHER DIRECTIONS}
\medskip

\medskip
{\bf 8.1. Large scale geometry.} 
The idea of studying properties of metric spaces at large scales was
introduced by Gromov [Gro81] in the context of groups of polynomial
growth. It was later developed into a broad framework for coarse geometry
and large scale geometry, see [NoYu12], [Roe03]. Certain (co)homology
functors for coarse geometry have been introduced in [Roe93], see also
Chapter~7 of [NoYu12]. In particular, the coaresening of homology theories 
described in Section~7.5 of [NoYu12] is based on the same notion of scale-dependent
Vietoris--Rips complexes that we discussed above in the setting of persistent
homology. Thus, we expect that the approach to persistence in terms of
Nori diagrams that we advocate in this paper should be applicable also to the context
of coarse geometry. It would be interesting to compare it with the axiomatic formulations
of coarse homology given in [Mitch01]. Among the interesting current problems in coarse
geometry are various topological and geometric rigidity conjectures (see Chapter~8
of [NoYu12]), which can be approached via index theory methods, developed in
the coarse geometry setting in  [Roe93]. It would be interesting to investigate
whether index theory in the coarse geometry context can be formulated in terms
of a more ``motivic" view of large scale geometry and coarse homology.

\medskip
{\bf 8.2. Cantor-like barcodes and fractality.}
It was obsevred in [KapVa04], [Pre11] how ind-pro objects
over a category behave as Cantor-like objects. This property
was used in [Li11] to model algebro-geometrically, in terms of
ind-pro varieties, the  energy-crystal momentum dispersion relation
for Harper and almost Mathieu operators with irrational parameters, 
replacing the ordinary Bloch variety by an ind-pro object, which
parallels the occurrence of the Hofstadter butterfly at the level of 
the spectrum, with its Cantor set fractal structure.
The density of states and the spectral functions are obtained in [Li11] 
as periods on this ``fractal-like" ind-pro version of the Bloch variety.

\smallskip

Within the context of this paper, one can consider the possibility of
extending the persistence structures and associated barcode diagrams
to a larger class of objects obtained as limits of finite type objects in
Vec$^{(\bold{R},\leq)}$ taken in such a way that the associated barcode
diagrams become Cantor sets. 
It would be interesting to investigate whether a larger class of physical
models similar to the algebro-geometric formulation of Harper operators 
given in [Li11] could be analyzed in terms of such limits of persistent 
homologies.

\medskip
{\bf 8.3. Khovanov homology and thin poset (co)homologies.}
In [Khov00] Khovanov constructed a categorification of the Jones polynomial
in the form of a chain complex of graded vector spaces and corresponding
homology whose graded Euler characteristic is the Jones polynomial. 
The chain complex is constructed by assigning to a knot or link
diagram with $N$ cossings the poset given by the $N$-cube and
a functor from this thin category to the category of graded vector 
spaces. The graded vector space assigned to a vertex of the cube
corresponds to a smoothing of the link where all the crossings are 
eliminated resulting in a union of $k$ planar closed curves, and the
graded dimension of the associated vector space depends on $k$ and
on the degree of the vertex, see [BarNat02] for more details. 

\smallskip

A generalization of this construction is given in [Chand18], where
a chain complex and a Khovanov-type cohomology $H^*(F,S,\Cal{A})$ 
are associated to any functor $F: (S,\leq) \to \Cal{A}$ from a poset to 
an abelian category, where the poset has a ``thinness" property 
described as follows. One requires the existence of a grading $r:S\to\bold{N}$
with $r(x)\leq r(y)$ for $x\leq y$, such that any pair $x,y\in S$ with $x\leq y$ for
which there is no $z\in S$ with $x< z <y$ should have $r(y)=r(x)+1$ and when
$r(y)=r(x)+2$ the set $\{z\in S\,:\, x< z <y \}$ consists of exactly two elements.
For posets $(S,\leq)$ that satisfy this thinness property, given a functor
$F: (S,\leq) \to \Cal{A}$ one constructs a chain complex with
$C^k(F,S,\Cal{A})=\oplus_{r(x)=k} F(x)$ and $\delta^k =
\sum c(x,y) F(x\leq y)$, where the sum is over all pairs $x,y$ with $x\leq y$ 
such that there is no $z$ with $x< z <y$, and $c$ is a ``balanced coloring".
This is a $\{ \pm 1 \}$ valued function on the set of pairs as above, with
the property that it has an odd number of $-1$'s on each ``diamond" set
$\{z\in S\,:\, x\leq z \leq y \}$ with $r(y)=r(x)+2$. The thinness property of
the poset and the balanced coloring property ensure that $\delta^2=0$
so that one obtains a chain complex. 
This general formalism is very suitable for introducing persistent versions of
Khovanov homology and related constructions and investigating the
topological information about knots and links that these persistent 
functors would capture.

\bigskip

\noindent {\bf Acknowledgment.} We thank Jack Morava for suggesting the question
of model structures for persistent homology discussed in Section~6. 
The second author is partially supported by NSF grant DMS-1707882,  
by NSERC Discovery Grant RGPIN-2018-04937 and Accelerator
Supplement Grant RGPAS-2018-522593, by the FQXi grant FQXi-RFP-1 804, 
and by the Perimeter Institute for Theoretical Physics.

\bigskip

\centerline{\bf References}

\medskip

[An04] Y.~Andr\'e. {\it Une introduction aux motives (motifs purs, motifs mixtes,
p\'eriodes.)} Panoramas et Synth\`eses, vol.~17. Soci\'et\'e Math\'ematique de France,
Paris, 2014.

\smallskip

[Ar13] D.~Arapura. {\it An abelian category of motivic sheaves.} Adv.~Math., 233, 2013,
pp.~135-195. arXiv:0801.0261

\smallskip

[Bar94]  S.A.~Barannikov, {\it The Framed Morse complex and its invariants}, Adv. Soviet Math. Vol.21 (1994)
93--115.

\smallskip

[BarNat02] D.~Bar-Natan, {\it On Khovanov's categorification of the Jones polynomial}, 
Algebraic \& Geometric Topology, 2 (2002) 337--370.  arXiv:math.QA/0201043

\smallskip

[BGSV90]  A.~Beilinson, A.~Goncharov, V.~Schechtman, A.~Varchenko. {\it
Aomoto dilogarithms, mixed Hodge structures and motivic cohomology of pairs of triangles 
on the plane}. In ``The Grothendieck Festschrift", Vol.~I, pp.135-172, Progr. Math., Vol.~86, 
Birkh\"auser, 1990.

\smallskip

[Be18] J.~Bergner, {\it The homotopy theory of $(\infty, 1)$-categories}. Cambridge University Press, 2018.

\smallskip
[BluLes17] A.~Blumberg, M.~Lesnick, {\it  Universality of the homotopy
interleaving distance}, arXiv:1705.01690

\smallskip

[BoChYv18] J.~Boissonnat, F.~Chazal, M.~Yvinec. {\it Geometric and topological inference}.
Cambridge University Press, 2018.

\smallskip

[Br14a] M.~Brion. {\it On algebraic semigroups and monoids}. In ``Algebraic monoids, group embeddings, and algebraic combinatorics", pp.1--54, Fields Inst. Commun., 71, Springer, 2014.

\smallskip

[Br14b] M.~Brion. {\it  On algebraic semigroups and monoids, II}. 
Semigroup Forum 88 (2014) no. 1, pp.~250-272.

\smallskip

[BuSc14] P.~Bubenik, J.~Scott. {\it Categorification of persistent homology.}
Discrete and Computational Geometry, 51 (3), 2014, pp.~600-627.

\smallskip

[BuSiSc15] P.~Bubenik, V.~de Silva, J.~Scott. {\it Metrics for generalised persistence
modules.}  Foundations of  Computational Math., vol.~15 (2015), issue 6, pp.~1501-1531. 
\smallskip

[Car09] G.~Carlsson, {\it  Topology and data}. 
Bull. Amer. Math. Soc. (N.S.) 46, no. 2, 2009, pp.~ 255-308. 

\smallskip

[Chand18] A.~Chandler, {\it Thin posets and homology theories}, preprint, 2018

https://alexchandler.wordpress.ncsu.edu/preprints/

\smallskip
[Du01] D.~Dugger, {\it Universal homotopy theories}. Adv. Math., Vol.164 (2001), N.1, 
pp.~144-176.
\smallskip
[EdHar10] H.~Edelsbrunner, J.~Harer. {\it  Computational topology}. American
Mathematical Society, 2010.
\smallskip
[EPY03] D.~Eisenbud, S.~Popescu, S.~Yuzvinsky. {\it Hyperplane arrangement cohomology and monomials in the exterior algebra}. Trans. Amer. Math. Soc. 355 (2003), no. 11, 4365-4383.

\smallskip

[GeMa03] S.~Gelfand, Yu.~Manin. {\it Methods of homological algebra.} 2nd Edition. Springer 2003, xvii + 372 pp.
\smallskip

[GoJa] P.~Goerss, R.~Jardine. {\it Simplicial homotopy theory.}  Birkh\"auser 1999.

\smallskip

[Gro81] M.~Gromov, {\it Groups of polynomial growth and expanding maps}, 
Publ. Math. IHES, 53 (1981) 53--73.

\smallskip

[Hat02] A.~Hatcher. {\it Algebraic Topology.} Cambridge UP, 2002.

\smallskip
[Hi03]  Ph.~Hirschhorn. {\it Model categories and their localizations.} American Mathematical Society,
2003.

\smallskip
[Ho98] M. Hovey, Model categories, Mathematical Surveys and Monographs, Vol. 63, American Math-
ematical Society, 1998.
\smallskip

[HuM-S17] A.~Huber, St.~M\"uller--Stach. {\it Periods and Nori motives. With
contributions by Benjamin Friedrich and Jonas von Wangenheim.} Springer, 2017,
xxiii+372 pp.

\smallskip

[KapVas04] M.~Kapranov, E.~Vasserot, {\it Vertex algebras and the formal loop space}, 
Publ. Math. IHES 100 (2004) 209--269.

\smallskip

[KaSch06] M.~Kashiwara, P.~Schapira. {\it Categories and sheaves.} Springer 2006, x + 497 pp.

\smallskip

[KaSch17] M.~Kashiwara, P.~Schapira. {\it Persistent homology and microlocal
sheaf theory.} arXiv:1705.00955 , 30 pp.

\smallskip

[Khov00] M.~Khovanov, {\it A categorification of the Jones polynomial}, Duke Math. J. 
101 (2000) no.~3, 359--426.  arXiv:math.QA/9908171

\smallskip

[Li11] D.~Li, {\it The algebraic geometry of Harper operators}, J. Phys. A 44 (2011), no. 40, 405204, 27 pp.

\smallskip

[Mar18] M. Marcolli, {\it Gamma Spaces and Information.} 
 J. Geom. Phys. 140 (2019), 26--55. 
arXiv:1807.05314.
 
\smallskip

[MaBo07] Yu.~Manin, D.~Borisov. {\it Generalized operads and their inner cohomomorhisms .}
 In: Geometry and Dynamics of Groups
and spaces (In memory of Aleksader Reznikov). Ed. by M. Kapranov et al.
Progress in Math., vol. 265. 
Birkh\"auser, Boston, 2007, pp.~247-308.
arXiv:math.CT/0609748

\smallskip

[Mitch01] P.D.~Mitchener, {\it Coarse homology theories}. Algebr. Geom. Topol. 1 (2001) 271--297.

\smallskip

[Mi86] H.~Mitsch. {\it A natural partial order for semigroups.} Proc. American Math. Soc., Vol. 97 (1986)
No. 3, 384-388.
\smallskip

[Mu06] D.~Murfet. {\it Abelian Categories.} Preprint, 2006

http://therisingsea.org/notes/AbelianCategories.pdf

\smallskip
[Na80]  K.~Nambooripad. {\it The natural partial order on a regular semigroup.} Proc. Edinburgh Math.
Soc. 23 (1980), pp.~ 249-260.

\smallskip

[NoYu12] P.W.~Nowak, G.~Yu, 
{\it Large scale geometry}, 
European Mathematical Society, 2012. xiv+189 pp.

\smallskip

[Pre11] L.~Previdi, {\it Locally compact objects in exact categories}, 
Internat. J. Math. 22 (2011) no. 12, 1787--1821.

\smallskip

[Roe03] J.~Roe, {\it Lectures on coarse geometry}, University Lecture Series, Vol.~31, 
American Mathematical Society, 2003. viii+175 pp.

\smallskip

[Roe93] J.~Roe, {\it Coarse cohomology and index theory on complete Riemannian manifolds}. 
Mem. Amer. Math. Soc. 104 (1993) no.~497, x+90 pp.

\smallskip
[SchSh00] S.~Schwede, B.~Shipley. {\it Algebras and modules in monoidal model categories.} Proc. London
Math. Soc. (3) 80 (2000), no. 2, pp.~491-511. 

arXiv:math/9801082

\smallskip

[Se74] G.~Segal. {\it  Categories and cohomology theories.} Topology, Vol.13 (1974) 293-312.

\smallskip

[Za05] A.~Zomorodian. {\it Topology for computing.} Cambridge University Press, 2005.
 
\enddocument